\documentclass[11pt,letterpaper,upref]{amsart}
\usepackage[letterpaper,margin=1.5in]{geometry}
\usepackage{amssymb}
\usepackage{amsxtra}
\usepackage[mathscr]{eucal}
\usepackage{graphics}
\usepackage{mathtools}
\usepackage[nobysame]{amsrefs}
\usepackage{color}
\usepackage{bm}
\makeatletter
\@namedef{subjclassname@2020}{\textup{2020} Mathematics Subject Classification}
\makeatother

\usepackage{caption}
\usepackage{subcaption}
\usepackage{tikz}
\usetikzlibrary{arrows.meta} 

\allowdisplaybreaks

\numberwithin{equation}{section}

\newtheorem{theorem}{Theorem}[section]
\newtheorem{proposition}[theorem]{Proposition}
\newtheorem{lemma}[theorem]{Lemma}
\newtheorem{corollary}[theorem]{Corollary}

\theoremstyle{definition}

\newtheorem{definition}[theorem]{Definition}
\newtheorem{example}[theorem]{Example}

\newtheorem{remark}[theorem]{Remark}
\newtheorem{remarks}[theorem]{Remarks}
\newtheorem{question}[theorem]{Question}

\renewcommand{\ge}{\geqslant}
\renewcommand{\le}{\leqslant}
\newcommand{\<}{\langle}
\renewcommand{\>}{\rangle}  

\renewcommand\Re{\operatorname{Re}}

\DeclareMathAlphabet{\mathsfit}{T1}{\sfdefault}{\mddefault}{\sldefault}
\SetMathAlphabet{\mathsfit}{bold}{T1}{\sfdefault}{\bfdefault}{\sldefault}

\newcommand{\af}{\alpha_f}
\newcommand{\Af}{\SA_f}
\newcommand{\Afo}{\SA_f^{\circ}}
\newcommand{\al}{\alpha} 
\newcommand{\SA}{\mathcal{A}}
\newcommand{\bo}{\boldsymbol{\omega}}
\newcommand{\bzero}{\boldsymbol{0}}
\newcommand{\bk}{\mathbf{k}}
\newcommand{\bfm}{\mathbf{m}}
\newcommand{\bn}{\mathbf{n}}
\newcommand{\br}{\mathbf{r}}
\newcommand{\bs}{\mathbf{s}}
\newcommand{\bu}{\mathbf{u}}
\newcommand{\bx}{\mathbf{x}}
\newcommand{\CC}{\mathbb{C}}
\newcommand{\czd}{\CC[\zd]}
\newcommand{\cnzd}{\CC[N\zd]}
\newcommand{\cst}{\CC^*}
\newcommand{\cstd}{(\cst)^d}
\newcommand{\ch}{\mathsfit{CH}}
\newcommand{\sfc}{^{\mathsf{c}}}

\newcommand{\Del}{\Delta}
\newcommand{\Delo}{\Delta^{\!\circ}}
\newcommand{\ssD}{\mathsfit{D}}
\newcommand{\DN}{\ssD_{N}}
\newcommand{\df}{\ssD_f}
\newcommand{\dlf}{\ssD_f^{}}
\newcommand{\dnf}{\ssD_N f}
\newcommand{\dom}{\operatorname{dom}}
\newcommand{\ssE}{\mathsfit{E}}
\newcommand{\EN}{\ssE_N}
\newcommand{\en}{e_N}
\newcommand{\epi}{\operatorname{epi}}
\newcommand{\eu}{e^{\bu}}
\newcommand{\et}{e^{2 \pi i \theta}}
\newcommand{\ep}{e^{2 \pi i \phi}}
\newcommand{\fa}{\mathfrak{a}} 
\newcommand{\fb}{\mathfrak{b}} 
\newcommand{\fc}{\mathfrak{c}}
\newcommand{\ck}{\fc_{\KK}}
\newcommand{\cl}{\fc_{\LL}}
\newcommand{\fp}{\mathfrak{p}}
\newcommand{\fq}{\mathfrak{q}}
\newcommand{\fh}{\widehat{f}}
\newcommand{\fnh}{\widehat{f}_{\<N\>}}
\newcommand{\fn}{f_{\<N\>}}
\newcommand{\fln}{f_{\<N\>}^{}}
\newcommand{\Gal}{\operatorname{Gal}}
\newcommand{\Gf}{\Gamma_{\! f}}
\newcommand{\Gp}{\Gamma_{\! \phi}}

\newcommand{\gh}{\widehat{g}}
\newcommand{\gn}{g_N}
\newcommand{\gnh}{\widehat{g}_N}
\newcommand{\gond}{G\ltimes \Omnd}
\newcommand{\gtn}{\widetilde{g}_N}
\newcommand{\h}{\mathsfit{h}} 
\newcommand{\Ht}{\mathsfit{H}}
\newcommand{\id}[1]{\< #1 \>}
\newcommand{\idr}[1]{\< #1 \>_R}
\newcommand{\ids}[1]{\< #1 \>_S}
\newcommand{\idzd}[1]{\< #1 \>_{\zzd}}
\newcommand{\idnzd}[1]{\< #1 \>_{\znzd}}
\newcommand{\KK}{\mathbb{K}}
\newcommand{\kzd}{\KK[\zd]}
\newcommand{\lzd}{\LL[\zd]}
\newcommand{\LL}{\mathbb{L}}
\newcommand{\ssL}{\mathsfit{L}}
\newcommand{\lam}{\lambda}
\newcommand{\Lam}{\Lambda}
\newcommand{\LN}{\ssL_{N}}
\newcommand{\mahler}{\mathsfit{m}}
\newcommand{\Mahler}{\mathsfit{M}}

\newcommand{\SN}{\mathcal{N}}
\newcommand{\nf}{\SN_f}
\newcommand{\nzd}{N\zd}
\newcommand{\OK}{\mathcal{O}_{\KK}}

\newcommand{\okzd}{\mathcal{O}_{\KK}[\zd]}

\newcommand{\OL}{\mathcal{O}_{\LL}}
\newcommand{\om}{\omega}

\newcommand{\Om}{\Omega}
\newcommand{\Omn}{\Om_N}
\newcommand{\Omnd}{\Om_N^d}
\newcommand{\omndp}{\Omnd\times\Phi}
\newcommand{\ph}{\widehat{\phi}}
\newcommand{\phir}{\phi_1\cdots\phi_r}

\newcommand{\QQ}{\mathbb{Q}}
\newcommand{\oR}{\overline{\RR}}
\newcommand{\ssR}{\mathsfit{R}}
\newcommand{\rd}{\RR^d}
\newcommand{\rf}{\ssR_f}
\newcommand{\rfs}{\ssR_f^*}
\newcommand{\RR}{\mathbb R}
\newcommand{\uR}{\underline{\RR}}
\newcommand{\SC}{\mathcal{C}}

\newcommand{\sn}{S_N}
\newcommand{\supp}{\operatorname{supp}}
\newcommand{\trop}{\operatorname{trop}}
\newcommand{\TT}{\mathbb{T}}
\newcommand{\tzd}{\TT^{\zd}}
\newcommand{\tnzd}{\TT^{N \zd}}
\newcommand{\USC}{\operatorname{USC}}
\newcommand{\variety}{\mathsfit{V}}
\newcommand{\xf}{X_f}
\newcommand{\ZZ}{\mathbb{Z}} 
\newcommand{\bz}{\mathbf{z}}
\newcommand{\zd}{\ZZ^{d}}
\newcommand{\zn}{\zeta_N}
\newcommand{\zxd}{\ZZ[x_1^{\pm1},\dots,x_d^{\pm 1}]}
\newcommand{\zzd}{\ZZ[\ZZ^d]}
\newcommand{\znzd}{\ZZ[N\zd]}

\begin{document}

\title{Decimation limits of principal algebraic $\zd$-actions}

\author[Arzhakova]{Elizaveta Arzhakova}
\address{Elizaveta Arzhakova: Mathematical Institute, Leiden University,
Postbus 9512, 2300 RA Leiden, The  Netherlands}
\email{e.arzhakova@math.leidenuniv.nl}

\author[Lind]{Douglas Lind}
\address{Douglas Lind: Department of Mathematics, University of
  Washington, Seattle, Washington 98195, USA}
  \email{lind@math.washington.edu}

\author[Schmidt]{\\Klaus Schmidt} 
\address{Klaus Schmidt: Mathematics Institute,
University of Vienna, Oskar-Morgenstern-Platz 1, A-1090 Vienna, Austria}
\email{klaus.schmidt@univie.ac.at}

\author[Verbitskiy]{Evgeny Verbitskiy}
\address{Evgeny Verbitskiy: Mathematical Institute, Leiden University,
Postbus 9512, 2300 RA Leiden, The  Netherlands  \newline\indent \textup{and}
\newline\indent Bernoulli Institute, University of Groningen, PO
Box 407, 9700 AK, Groningen, The Netherlands}
\email{evgeny@math.leidenuniv.nl}

\date{\today}




\begin{abstract}
   Let $f$ be a Laurent polynomial in $d$ commuting variables with integer
   coefficients. Associated to $f$ is the principal algebraic $\zd$-action $\af$ on a compact
   subgroup $X_f$ of $\TT^{\zd}$ determined by $f$. Let $N\ge1$ and restrict points in $X_f$ to
   coordinates in $N\zd$. The resulting algebraic $N\zd$-action is again principal, and is
   associated to a polynomial $g_N$ whose support grows with $N$ and whose coefficients grow
   exponentially with $N$. We prove that by suitably renormalizing these decimations we can
   identify a limiting behavior given by a continuous concave function on the Newton polytope
   of $f$, and show that this decimation limit is the negative of the Legendre dual of the
   Ronkin function of $f$. In certain cases with two variables, the decimation limit coincides with the surface tension of random surfaces related to dimer models, but the statistical physics methods used to prove this are quite different and depend on special properties of the polynomial.
\end{abstract}

\maketitle

\section{Introduction}\label{sec:introduction}

Let $d\ge1$ and $f\in\zxd$ be a Laurent polynomial with integer coefficients in $d$ commuting
variables. We write $f(x_1,\dots,x_d)=f(\bx)=\sum_{\bn\in\zd}\fh(\bn)\bx^{\bn}$, where
$\bx^{\bn}=x_1^{n_1}\dots x_d^{n_d}$ and $\fh(\bn)\in\ZZ$ for all $\bn\in\zd$ and is nonzero
for only finitely many $\bn\in\zd$.

Denote the additive torus $\RR/\ZZ$ by $\TT$. Use $f$ to define a compact subgroup $X_f$ of~
$\TT^{\zd}$ by
\begin{equation}\label{eqn:principal-action}
   X_f :=\Bigl\{ t\in\TT^{\zd}:\sum_{\bn\in\zd} \fh(\bn)t_{\bfm+\bn}=0 \text{\quad for all
   $\bfm\in\zd$} \Bigr\}.
\end{equation}
By its definition this subgroup is invariant under the natural shift-action $\sigma$ of $\zd$
on~ $\TT^{\zd}$ defined by $\sigma^{\bn}(t)_{\bfm}=t_{\bfm-\bn}$. Hence the restriction $\af$ of $\sigma$ to $X_f$ gives an action of $\zd$ by automorphisms of the compact abelian group
$X_f$. We call $(X_f,\af)$ the \emph{principal algebraic $\zd$-action defined by $f$}.

Such $\zd$-actions serve as a rich class of examples and have been studied intensively. An
observation of Halmos \cite{Halmos} shows that $\af$ automatically preserves Haar measure
$\mu_f$ on $X_f$. It is known that the topological entropy of $\af$ coincides with its
measure-theoretic entropy with respect to $\mu_f$. For nonzero $f$ this common value was computed in \cite{LSW} to be the logarithmic Mahler measure of $f$, defined as
\begin{equation}\label{eqn:mahler-measure}
   \mahler(f) := \int_0^1 \dots \int_0^1 \log|f(e^{2\pi i s_1},\dots, e^{2\pi i s_d})| \,ds_1\dots ds_d
\end{equation}
(when $f=0$ the entropy is infinite).

It will be convenient to identify the Laurent polynomial ring $\zxd$ with the integral group
ring $\zzd$, where the monomial $\bx^{\bn}$ corresponds to $\bn\in\zd$. Thus $f\in\zxd$ is
identified with its coefficient function $\fh\colon\zd\to\ZZ$. When emphasizing the behavior of
coefficients we will always use the notation $\fh$. 

Principal algebraic $\zd$-actions are special cases of \emph{algebraic $\zd$-actions}, which we define as actions of $\zd$ by automorphisms of a compact abelian group. Modules over the integral group ring $\zzd$ arise naturally as Pontryagin duals of algebraic $\zd$-actions: an algebraic $\zd$-action on $X$ induces a $\zzd$-module structure on the discrete dual group $\widehat{X}$, and conversely if $M$ is a $\zzd$-module there is an induced algebraic $\zd$-action on its compact dual group $\widehat{M}$ (see \cite{SchmidtBook}*{Chap.\ II} for a detailed description). In \S\ref{sec:algebraic-decimations} we give a concrete account of this correspondence for the algebraic actions we study here. 

Fix a principal algebraic $\zd$-action $(X_f,\af)$. Let $N\ge1$ and
$r_N\colon\TT^{\zd}\to\TT^{N\zd}$ be the map restricting the coordinates of a point to only those in the sublattice $N\zd$. We call the image $r_N(X_f)$ the \emph{$N$th decimation} of
$X_f$, although this is considerably more brutal that the term's original meaning since only
every $N$th coordinate in each coordinate direction survives. Clearly $r_N(X_f)$ is again a compact abelian group, and it is
invariant under the natural shift action of $N\zd$ on $\TT^{N\zd}$.

Using commutative algebra applied to contracted ideals in integral extensions, we show in
\S\ref{sec:algebraic-decimations} that $r_N(X_f)$ is a principal algebraic $N\zd$-action with
some defining polynomial $g_N\in\znzd$. Typically both the support of $g_N$ grows with $N$ and its
coefficient function $\gnh$ grows exponentially in $N$. Our goal in this paper is to prove that with suitable renormalizations the concave hulls of the resulting functions converge uniformly on the Newton polytope of $f$ to a continuous decimation limit $\df$. Furthermore, $\df$ can be computed via Legendre duality using a well-studied object called the Ronkin function of $f$.

The analytical parts of our analysis apply to Laurent polynomials with complex
coefficients. For such an $f\in\czd$ we define its $N$th \emph{decimation} $\fn$ by
\begin{equation}\label{eqn:decimation}
   \fn(x_1,\dots,x_d):= \prod_{k_1=0}^{N-1}\dots \prod_{k_d=0}^{N-1} f(e^{2 \pi i k_1/N}x_1,
   \dots, e^{2 \pi i k_d/N}x_d).
\end{equation}
Since $\fn$ is unchanged after multiplying each of its variables by an arbitrary $N$th root of unity, it follows that it is a polynomial in the $N$th powers of the $x_i$, i.e., that $\fn\in\cnzd$. Decimations of polynomials have appeared in many contexts, including Purbhoo's approximations to shapes of complex amoebas \cite{Purbhoo}, Boyd's proof that the Mahler measure of a polynomial is continuous in its coefficients \cite{BoydUniform}, and dimer models in statistical physics \cite{KenyonOkounkovSheffield}.

For most $f\in\zzd$ the generator $g_N$ of the $N$th decimation of $X_f$ coincides with~
$\fn$. But under special circumstances characterized in \S\ref{sec:algebraic-decimations},
involving  the support of $f$ and the Galois properties of the coefficients of the
polynomials occurring in the factorization of $f$ over the algebraic closure of the rationals, it can happen that $\fn$ is a power $g_N^{e_N}$ of $\gn$. Our results give a way to explicitly compute the exponent $e_N$ and provide a way to determine the renormalized behavior of $g_N$ from that of $\fn$.

To give a simple example when $d=1$, let $f(x)=x^2 -2$. Then since $f$ is already in $\ZZ[2\ZZ]$ we have that $g_2(x)=f(x)$, while
$f_{\<2\>}(x)=f(x)f(-x)=f(x)^2=g_2(x)^2$. In Example \ref{exam:1-2x2} we show that the exponents for this case are $e_N= 1$ if $N$ is odd, while $e_N=2$ if $N$ is even.

For $f\in\czd$ let $\supp f=\{\bn\in\zd:\fh(\bn)\ne0\}$ denote its support. The \emph{Newton
polytope} $\nf$ of $f$ is the convex hull in $\rd$ of $\supp f$. It is well-known that if $f$ and $g$ are in $\czd$, then $\SN_{fg}$ is the Minkowski sum $\SN_f+\SN_g$. 
Since $\fn$ is the product of $N^d$
polynomials all of whose Newton polytopes are $\nf$, it follows 
that $\SN_{\fn}=N^d\nf$.

The \emph{Ronkin function} $\rf\colon\rd\to\RR$ of $0\ne f\in\czd$ is defined by
\begin{equation}\label{eqn:ronkin}
   \rf(u_1,\dots,u_d) := \int_0^1\dots \int_0^1 \log|f(e^{u_1}e^{2 \pi i s_1},\dots,
   e^{u_d} e^{2\pi i s_d})|\, ds_1\dots ds_d .
\end{equation}
By \cite{PassareRullgard} this a convex function on $\rd$. Letting $\<\br,\bu\>= r_1u_1+\dots +r_du_d$ denote the usual inner product on $\rd$, then $\rf$ has a Legendre dual $\rfs$ defined by
\begin{displaymath}
   \rfs(\br) := \sup \{\<\br,\bu\>-\rf(\bu): \bu\in\rd\},
\end{displaymath}
which turns out to be a convex function on $\nf$ (and is $\infty$ off $\nf$).

To describe rescaling of polynomials $g\in\czd$ it is convenient to extend the domain of $\gh$
from $\zd$ to $\rd$ by declaring its value to be 0 off $\supp g$. 

Let $\phi\colon\rd\to\CC$. For any $a>0$ define  the \emph{rescaling operator} $\ssE_a$ on $\phi$ by  $(\ssE_a\phi)(\br) = \phi(a\br)$ for all $\br \in\rd$. When dealing with concave functions it is often convenient to use the extended range $\uR=\RR\cup\{-\infty\}$, with the usual algebraic rules for handling $-\infty$ and with the convention that $\log 0 = -\infty$. Then $\log|\phi|\colon\rd\to \uR$, and we define its \emph{concave hull} $\ch(\log|\phi|)$ to be the infimum of all affine functions on $\rd$ that dominate $\log|\phi|$.

\begin{definition}\label{def:decimation-limit}
   Let $f\in\czd$ and $\fn$ be its $N$th decimation. Define the $N$th \emph{logarithmic rescaling}
   $\LN f$ of $f$ by
   \begin{displaymath}
      \LN f:= \ssE_{N^d} \Bigl( \frac{1}{N^d}\log |\fnh| \Bigr).
   \end{displaymath}
   Clearly $\LN f(\br)=-\infty $ if $\br\notin \nf$, and is finite at every extreme point of
   $\nf$ and at only finitely many other points in $\nf$. The $N$th \emph{renormalized
   decimation} $\DN f$ of $f$ is the concave hull $\ch(\LN f)$  of $\LN f$. By our previous remark, $\DN f$ equals $-\infty$ off $\nf$ and is finite at every point of $\nf$.
\end{definition}

With these preparations we can now state one of our main results.

\begin{theorem}\label{thm:main}
   Let $0\ne f \in \czd$. Then the $N\!$th renormalized decimations $\DN f$ of $f$ are concave polyhedral functions on the Newton polytope $\nf$ of $f$ that
   converge uniformly on $\nf$ as $N\to\infty$ to a continuous concave  limit function ~$\df$ called the \textbf{decimation limit} of~$f$. Off $\nf$ both $\DN f$ and $\df$ are equal to $-\infty$. Furthermore $\ssD_f^{}=-\rfs$, where $\rfs$ is the Legendre dual of the Ronkin function $\rf$ of $f$.
\end{theorem}

The proof of this theorem uses two main ideas: Mahler's fundamental estimate \cite{Mahler}
relating the largest coefficient of a polynomial to its Mahler measure and support, and a method used by Boyd \cite{BoydUniform}, applied to decimations along powers of 2, to prove that for polynomials whose support is contained in a fixed finite subset of $\zd$ the Mahler measure is a continuous function of their coefficients.

If $f\in\zzd$ the decimation limit of $f$ contains dynamical information about ~ $\af$.

\begin{corollary}\label{cor:max}
   Let $0\ne f\in\czd$. Then the maximum value of the decimation limit $\df$ on the Newton polytope $\nf$ equals the logarithmic Mahler measure $\mahler(f)$ of ~$f$ defined in \eqref{eqn:mahler-measure}. In particular, if $f\in\zzd$ then this maximum value equals the entropy of the principal algebraic $\zd$-action $\af$.
\end{corollary}

Duality allows us to compute the decimation limit of a product of two polynomials. Suppose that
$\phi,\psi\colon\rd\to\uR$ both have finite supremum. Define their \emph{tropical convolution}
$\phi\circledast\psi$ by
\begin{displaymath}
   (\phi\circledast\psi)(\br) := \sup\{\phi(\bs)+\psi(\br-\bs):\bs\in\rd \}.
\end{displaymath}
This is the tropical analogue of standard convolution, but using tropical (or max-plus)
arithmetic in $\uR$.

\begin{corollary}\label{cor:tropical-convolution}
   Let $f$ and $g$ be nonzero polynomials in $\czd$. Then $\ssD_{fg}=\ssD_f \circledast \ssD_g$.
\end{corollary}

Thus decimation limits live in the tropics.

\medskip
The authors are grateful to Hanfeng Li for several suggestions and clarifications.

\section{Examples}\label{sec:examples}

Here we give some examples to illustrate the  phenomena we are investigating. They use either one or two variables, and for these we denote the variables by $x$ and ~$y$ rather than
$x_1$ and $x_2$. Let $\Omn=\{e^{2\pi i k/N}:0\le k < N\}$ denote the group of $N$th roots of unity.

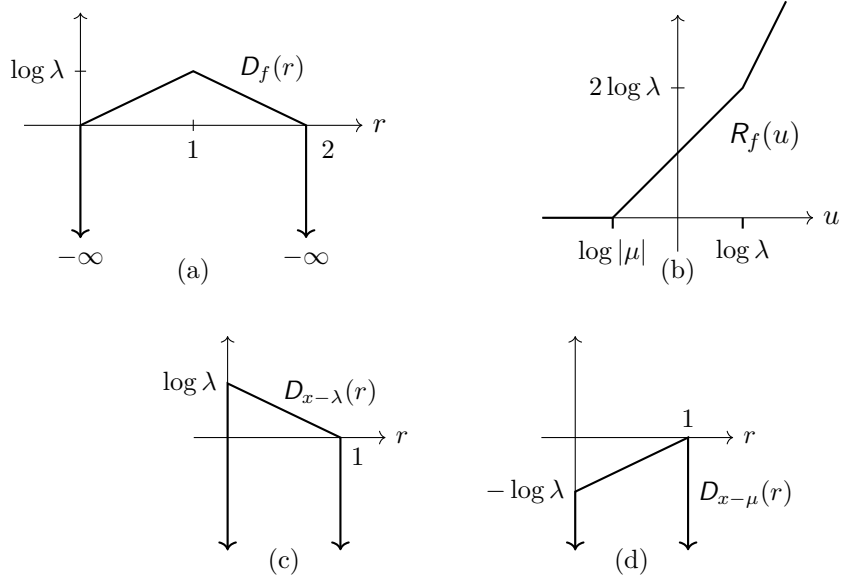
\begin{figure}[ht]
   \centering
   \begin{subfigure}[b]{0.4\textwidth}
   	 \centering
       \begin{tikzpicture}[scale=1.5]
       	   \draw[->,thin] (0,-1) -- (0,1);
       	   \draw[->,thin] (-0.5,0) -- (2.5,0) node[right] {$r$};
       	   \draw[line width=0.8pt,<->] (0,-1) -- (0,0) -- (1,.48) -- (2,0) -- (2,-1);
       	   \draw (0,-1) node[below]  {\small$-\infty$};
       	   \draw (2,-1) node[below]  {\small$-\infty$};
       	   \draw (1,-0.05) node[below]{\small1}-- (1,0.05);
       	   \draw (-0.05,.48) node[left] {\small$\log \lambda$} -- (0.05,.48);
       	   \draw (2,0) node[below right=2pt]{\small 2};
       	   \draw (1.7,.5) node{\small$\ssD_f(r)$};
       	   \draw (1,-1.3) node{\small (a)};
       \end{tikzpicture}
   \end{subfigure}
   \hfil
	\begin{subfigure}[b]{0.4\textwidth}
	   \centering
	   \begin{tikzpicture}[scale=1.8]
	       \draw[->,thin] (0,-0.25) -- (0,1.5);
	       \draw[->,thin] (-1,0) -- (1,0) node[right] {$u$};
	       \draw[line width=0.8pt] (-1,0) -- (-0.48,0) -- (.48,.96) -- (.8,1.6);
	       \draw[thick] (-0.48,-0.08) node[below]{\small $\log |\mu|$} -- (-0.48,0);
	       \draw[thick] (0.48,-0.08) node[below]{\small $\log \lambda$} -- (0.48,0);
	       \draw (-.05,.96) node[left]{\small $2\log \lambda$} -- (0.05,0.96);
	       \draw (0.3,0.6) node[right]{$\rf(u)$};
	       \draw (0,-.4) node{\small (b)};
	   \end{tikzpicture}
	\end{subfigure}\vskip 0.2in
	\begin{subfigure}[b]{0.4\textwidth}
		\hfill
		\begin{tikzpicture}[scale=1.5]
		    \draw[->,thin] (0,-1) -- (0,0.9);
		    \draw[->,thin] (-0.3,0) -- (1.4,0) node[right] {\small $r$};
		    \draw[line width=0.8pt,<->] (0,-1) -- (0,0.48) node[left] {\small $\log \lambda$} -- (1,0) node[below right]{\small 1} -- (1,-1);
		    \draw (.4,.4) node[right]{\small$\ssD_{x-\lambda}(r)$};
		    \draw (0.5,-1.1) node{\small (c)};
		\end{tikzpicture}
	\end{subfigure}\quad\quad
	\begin{subfigure}[b]{0.4\textwidth}
		\begin{tikzpicture}[scale=1.5]
			 \draw[->,thin] (0,-1) -- (0,0.9);
			 \draw[->,thin] (-0.3,0) -- (1.4,0) node[right] {\small $r$};
			 \draw[line width=0.8pt,<->] (0,-1) -- (0,-0.48) node[left] {\small $-\log \lambda$} -- (1,0) node[above]{\small 1} -- (1,-1);
			 \draw (1,-.5) node[right]{\small$\ssD_{x-\mu}(r)$};
			 \draw (0.5,-1.1) node{\small (d)};
		\end{tikzpicture}
		\hfill
	\end{subfigure}
	\caption{Graphs in Example \ref{exam:golden}}\label{fig:golden}
\end{figure}

\begin{example}\label{exam:golden}
   Let $d=1$ and $f(x)=x^2-x-1=(x-\lambda)(x-\mu)$, where $\lambda=(1+\sqrt{5})/2$ and
   $\mu=(1-\sqrt{5})/2$. Then
   \begin{align*}
      \fn(x)&=\prod_{\om\in\Omn}f(\om x)=\prod_{\om\in\Omn}(\om x-\lambda)(\om x-\mu)\\
            &=(x^N-\lambda^N)(x^N-\mu^N)=x^{2N}-(\lambda^n+\mu^N)x^N + (-1)^N.
   \end{align*}
   Hence
   \begin{displaymath}
      (\LN f)(r) =
      \begin{cases}
         0       & \text{if $r=0$ or $2$}, \\
         \frac{1}{N} \log |\lambda^N + \mu^N| &\text{if $r=1$}, \\
         -\infty       &\text{otherwise} .
      \end{cases}
   \end{displaymath}
   Since $\LN f(1)\to \log\lambda$ as $N\to\infty$, the concave hulls $\DN f$ converge
   uniformly on $\nf=[0,2]$ to the decimation limit
   \begin{displaymath}
      \df(r) =
      \begin{cases}
         r \log\lambda  &\text{if $0\le r\le 1$}, \\
         (2-r)\log \lambda   &\text{if $1\le r\le 2$}, \\
         -\infty   &\text{otherwise},
      \end{cases}
   \end{displaymath}
   which is shown in Figure \ref{fig:golden}(a).
   
	To compute the Ronkin function $\rf$, recall Jensen's formula that for every $\xi\in\CC$ we have that
	\begin{equation}\label{eqn:jensen}
	   \int_0^1 \log |e^{2\pi i s}-\xi|\,ds = \max\{0,\log|\xi|\} :=\log^+|\xi|.
	\end{equation}
	Thus
	\begin{align*}
		\rf(u) = \int_0^1 \log|f(e^u e^{2\pi i s})|\,ds &= \int_0^1\log|e^u e^{2\pi i s}-\lambda|\,ds + \int_0^1 \log|e^u e^{2\pi i s}-\mu|\, ds \\
		&= 2u + \log^+|e^{-u}\lambda|+\log^+|e^{-u}\mu|,
	\end{align*}
	whose polygonal graph is depicted in Figure \ref{fig:golden}(b). It is then easy to verify using the definition of Legendre transform that $\ssD_f^{}=-\rfs$.
	
	Finally, the decimation limits $\ssD_{x-\lambda}$ and $\ssD_{x-\mu}$ are computed similarly, and shown in Figures \ref{fig:golden}(c) and \ref{fig:golden}(d). It is easy to check using the definition of tropical convolution that $\ssD_{x-\lambda}\circledast \ssD_{x-\mu} = \ssD_{(x-\lambda)(x-\mu)} = \df$, in agreement with Corollary \ref{cor:tropical-convolution}. 
\end{example}

More generally, if $f(x)=\prod_{j=1}^m (x-\lambda_j)$ and $|\lam_1|>|\lam_2|>\dots>|\lam_m|$, then a computation similar to that in Example \ref{exam:golden} shows that 
$(\LN f)(m)=0$ and that 
$(\LN f)(k)$ converges to $\log|\lam_1\lam_2 \dots\lam_{m-k}|$ for $k=0,1,\dots,m-1$, and this gives uniform convergence of $\dnf$ to $\df$ on $\nf=[0,m]$. However, if some roots of $f$ have equal absolute value, then convergence is more delicate, or may even fail, as the next two examples show.

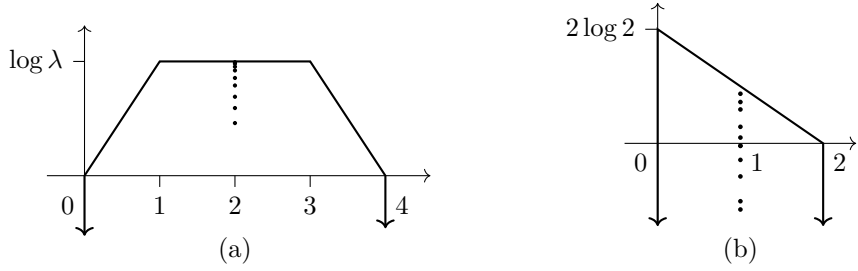
\begin{figure}[ht]
   \centering
   \begin{subfigure}[b]{0.4\textwidth}
   	 \centering
       \begin{tikzpicture}
       	   \draw[->] (0,-0.7) -- (0,2);
       	   \draw[->] (-0.5,0) -- (4.6,0);
       	   \draw[<->,line width=0.8pt] (0,-0.8) -- (0,0) --(1,1.52)--(3,1.52)--(4,0)-- (4,-0.7);
       	   \draw[] (0,1.52) -- (-0.15,1.52) node[left]{\small $\log\lambda$};
       	   \draw (0,-0.15) node[below left]{\small 0};
       	   \draw (1,0) -- (1,-0.15) node[below]{\small 1};
       	   \draw (2,0) -- (2,-0.15) node[below]{\small 2};
       	   \draw (3,0) -- (3,-0.15) node[below]{\small 3};
       	   \fill[black] (2,1.3) circle (0.8pt);
       	   \fill[black] (2,1.4) circle (0.8pt);
       	   \fill[black] (2,1.45) circle (0.8pt);
       	   \fill[black] (2,1.49) circle (0.8pt);
       	   \fill[black] (2,1.51) circle (0.8pt);     	   
       	   \fill[black] (2,0.7) circle (0.8pt);
       	   \fill[black] (2,0.9) circle (0.8pt);
       	   \fill[black] (2,1.2) circle (0.8pt);
       	   \fill[black] (2,1.05) circle (0.8pt);
       	   \draw (4,-0.15) node[below right]{\small 4};
       	   \draw (2,-1) node{\small (a)};
       \end{tikzpicture}
   \end{subfigure}
   \hfil
	\begin{subfigure}[b]{0.4\textwidth}
	   \centering
	   \begin{tikzpicture}[scale=1.1]
	       \draw[->] (-0.4,0) -- (2.4,0);
	       \draw[->] (0,-1) -- (0,1.7);
	       \draw[<->, line width = 0.8pt] (0,-1)--(0,1.38)--(2,0) -- (2,-1);
	       \draw[] (0,1.38) -- (-0.15,1.38) node[left]{\small $2 \log 2$};
			 \fill[black] (1,0.6) circle (0.8pt); 
			 \fill[black] (1,0.5) circle (0.8pt);
			 \fill[black] (1,0.41) circle (0.8pt);
			 \fill[black] (1,0.2) circle (0.8pt);
			 \fill[black] (1,0.07) circle (0.8pt);
			 \fill[black] (1,-0.03) circle (0.8pt);
			 \fill[black] (1,-.4) circle (0.8pt);
			 \fill[black] (1,-.7) circle (0.8pt); \fill[black] (1,-0.03) circle (0.8pt);
			 \fill[black] (1,-.8) circle (0.8pt);
			 \fill[black] (1,-.2) circle (0.8pt); 
			 \draw (2,0) node[below right]{\small 2};
			 \draw (0,0) node[below left]{\small 0};
			 \draw (1,0) node[below right]{\small 1};
			 \draw (1,-1.3) node{\small (b)};
	   \end{tikzpicture}
	\end{subfigure}
	\caption{(a) Convergence in Example \ref{exam:quasihyperbolic}, and (b) lack of convergence
	   in Example \ref{exam:transcendental} \label{fig:nonconvergence}}
\end{figure}

\begin{example}\label{exam:quasihyperbolic}
   Let $d=1$ and $f(x)=x^4-4x^3-2x^2-4x+1$, which is irreducible in $\ZZ[\ZZ]$. The roots of
   $f$ are $\lam = 1+\sqrt{2} + \sqrt{2\sqrt{2} + 2}\approx 4.611$,  $\mu=1+\sqrt{2} -
   \sqrt{2\sqrt{2}+2} \approx 0.217$, and $1-\sqrt{2}\pm i \sqrt{2\sqrt{2}-2} = e^{\pm 2\pi i
   \theta}$, where $\theta$ is irrational. Simple estimates show that $(\LN f)(k)$ converges
   for $k=0,1,3,4$ with limits $0,\log\lam,\log\lam,0$, respectively. However, the dominant
   term controlling the behavior of $(\LN f)(2)$ is
   \begin{displaymath}
      \frac{1}{N}\log|2\lam^N \cos(2\pi N \theta)|.
   \end{displaymath}
   Since $\theta$ is irrational, the factor $\cos(2\pi N\theta)$ occasionally becomes very
   small, and so convergence is in question.
   
   In fact, $(\LN f)(2)$ does converge, but the proof requires a deep result of Gelfond~ \cite{Gelfond}*{Thm.\ III, p.\ 28} on the diophantine properties of algebraic numbers on the unit circle (see \cite{Baker}*{Thm.\ 3.1} for a more accessible treatment). According to this result, if $\xi$ is an algebraic number (such as $e^{2\pi i \theta}$ above) such that $|\xi|=1$ and $\xi$ is not a root of unity, and if $\varepsilon >0$, then $|\xi^n-1|>e^{-n\varepsilon}$ for all but finitely many ~$n$. From this it is easy to deduce that $|e^{2\pi i N \theta}-i|>e^{-N\varepsilon}$ for almost every $N$, and hence that $(1/N)\log|\cos(2\pi N\theta)|\to0$ as $N\to \infty$. This convergence is illustrated in Figure \ref{fig:nonconvergence}(a).

   Both $(\LN f)(1)$ and $(\LN f)(3)$ converge to $\log \lam$, and clearly
   $\limsup_{N\to\infty} (\LN f)(2)\le\log\lam$. Hence any lack of convergence of $(\LN f)(2)$
   would not affect the limiting behavior of the concave hull $\dnf$, nor uniform convergence of $\dnf$ to $\df$ on $[0,4]$. Thus such diophantine issues are covered up by taking concave hulls.
   
\end{example}

The next example shows that if we allow the coefficients of $f$ to be arbitrary complex numbers instead of integers, then $(\LN f)(k)$ can badly fail to converge at some $k$. 

\begin{example} \label{exam:transcendental}
   Let $d=1$ and $f(x)=(x-2e^{2\pi i \theta})(x - 2e^{-2\pi i \theta})$, where we will
   determine $\theta$. Then $(\LN f)(0)=2\log 2$ and $(\LN f)(2)= 0$ for all $N\ge 1$, while
   \begin{displaymath}
      (\LN f)(1) = \frac{1}{N} \log | 2^N \cdot 2\cos(2\pi N \theta)|.
   \end{displaymath}
   It is possible to construct an irrational  $\theta$ and a sequence $N_j\to\infty$
   such that $\frac{1}{N_j} \log|\cos(2\pi N_j \theta)|\to-\infty$ as $j\to\infty$. Hence
   using this value of $\theta$ to define $f$ we see that $(\LN f)(1)$ does not converge, as depicted in Figure \ref{fig:nonconvergence}(b), although the concave hulls $\dnf$ do converge uniformly to $\df$.
\end{example}

Using arguments similar to those above, it is possible to give an elementary direct proof of Theorem \ref{thm:main} in the case $d=1$.

\begin{example}\label{exam:ledrappier}
   Let $d=2$ and $f(x,y)=1+x+y$. Then $\fn$ is a polynomial in $x^N$ and $y^N$ of degree
   $N^2$ in each variable. For example,
   \begin{align*}
      f_{\<5\>}(x,y) =\,\,  & x^{25}  +5 x^{20} y^5+5 x^{20}+10 x^{15} y^{10}-605 x^{15} y^5+10 x^{15}+10
      x^{10} 
      y^{15} \\
      & \quad +1905 x^{10} y^{10}+1905 x^{10}
      y^5+10 x^{10}+5 x^5 y^{20} -605
      x^5 y^{15}  + 1905 x^5 y^{10}\\
      &  \quad\quad  -
      605 x^5 y^5+5 x^5+y^{25}+5 y^{20}+10 y^{15}+10 y^{10}+5 y^5+1.
   \end{align*}
   The $N$th logarithmic rescaling $\LN f$ of $f$ is finite at points in the unit simplex
   $\Delta=\nf$ whose coordinates are integer multiples of $1/N$. Thus its concave hull
   $\dnf$ is a polyhedral surface over $\Delta$, and as $N\to\infty$ these surfaces converge
   uniformly on $\Delta$ to the graph of the concave decimation limit $\df$. Figure \ref{fig:ledrappier}(a) shows the polyhedral surface $D_5 f$ corresponding to the calculation of $f_{\<5\>}$ above, and Figure~3(b) depicts the limiting smooth surface for $\df$.
   
   	\begin{figure}[b]
      	\begin{subfigure}[b]{0.70\textwidth}
      		\centering
         		\includegraphics[width=\textwidth]{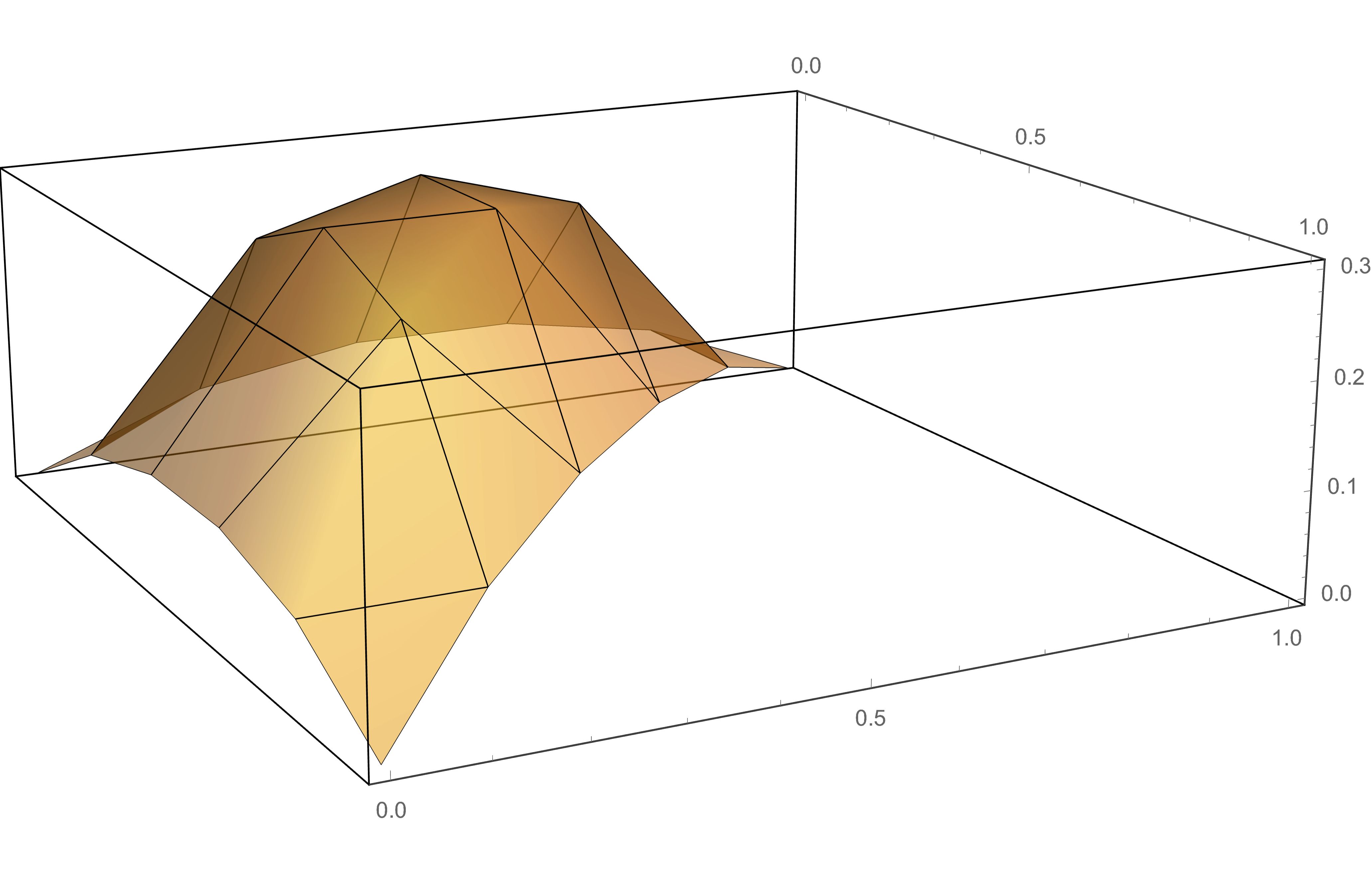}\caption{}
      	\end{subfigure}\\
   		\begin{subfigure}[b]{0.70\textwidth}
   			\centering
   	   	\includegraphics[width=1.0\textwidth]{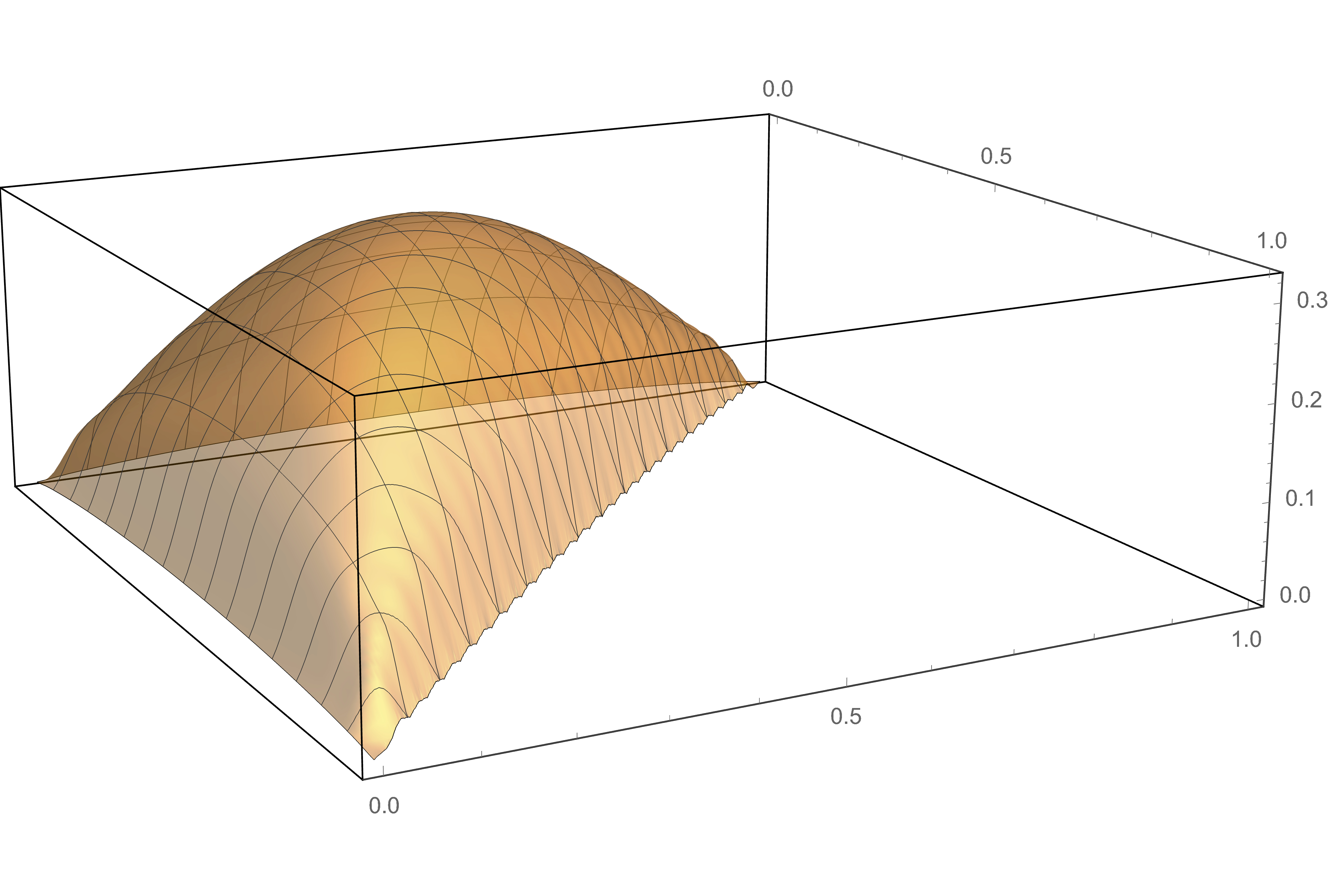}\caption{}
   		\end{subfigure}
   		\caption{(a) Polyhedral approximation $\ssD_5f$, and (b) limiting smooth surface $\df$ for $f(x,y)=1+x+y$ in Example \ref{exam:ledrappier}} \label{fig:ledrappier}
   	\end{figure}

   For this example it is possible to derive an explicit formula for $\df$. Clearly $\df(r,s)$
   is symmetric in $r$ and $s$, so we may assume that $s\le r$. Let
   \begin{align}
     \Delta_1 &= \{ (r,s)\in \Delta: s\le r \text{ and } s \le (1-r)/2\}, \label{eqn:delta1}\\
     \Delta_2 &= \{ (r,s)\in \Delta: s\le r \text{ and } s \ge (1-r)/2\}. \label{eqn:delta2}
   \end{align}
   For $(r,s)\in\Delta_1\cup\Delta_2$ with $r+s<1$ define
   \begin{displaymath}
   b(r,s)=\csc[\pi(r+s)]\sin(\pi s).
   \end{displaymath}
   Then it turns out that $0\le b(r,s)\le1$ for $(r,s)\in\Delta_1$ while $1\le b(r,s)<\infty$ for
   $(r,s)\in\Delta_2$.

   Using Legendre duality and calculations of $\rf$ by Lundqvist \cite{Lundqvist}, 
   in Appendix A we show that that if $(r,s)\in\Delta_1$ then
   \begin{equation}\label{eqn:led1}
      \df(r,s) =\sum_{n=1}^\infty \frac{(-1)^{n+1}}{\pi n^2} \, b(r,s)^n\sin[n \pi(1-r)]-s\log b(r,s),
   \end{equation}
   while if $(r,s)\in\Delta_2$ then
   \begin{equation}\label{eqn:led2}
      \df(r,s) = \sum_{n=1}^\infty \frac{(-1)^{n+1}}{\pi n^2} \, b(r,s)^{-n}\sin[n \pi(1-r)]+
      (1-r-s)\log b(r,s).
   \end{equation}

   We will prove in Corollary \ref{cor:max} that the maximum value of $\df$ equals the entropy of ~$\af$, which is the logarithmic Mahler measure $\mahler(f)$ of $f$ defined in \eqref{eqn:mahler-measure}. In this example, the maximum value is attained at $(1/3,1/3)$, which is in both $\Delta_1$ and $\Delta_2$. Either formula therefore applies, and each gives Smyth's calculation \cite{Smyth} that 
   \begin{equation}\label{eqn:smyth}
      \mahler(1+x+y)= \df(1/3,1/3) =\frac{3\sqrt{3}}{4\pi}\sum_{n=1}^\infty \frac{\chi_3(n)}{n^2}
      = \frac{3\sqrt{3}}{4\pi}L(2,\chi_3) \approx 0.3230,
   \end{equation}
   where $\chi_3$ is the nontrivial character of $\ZZ/3\ZZ$ and $L(s,\chi_3)$ is the $L$-function
   associated with $\chi_3$.
\end{example}

Unlike the previous example, some decimation limits exhibit non-smooth behavior.

\begin{example}\label{exam:amoeba-with-hole}
   Let $d=2$ and $f(x,y)= 5 + x + x^{-1}+y+y^{-1}$. The decimation limit $\df$ is depicted in
   Figure \ref{fig:amoeba-with-hole}(a). The non-smooth peak at the origin is due to
   a ``hole'' in the amoeba of $f$, as defined in \S\ref{sec:ronkin-functions} and shown in Figure~\ref{fig:amoeba-with-hole}(b).

	\begin{figure}[ht]
   	\begin{subfigure}[b]{0.5\textwidth}
   		\centering
      		\includegraphics[width=\textwidth]{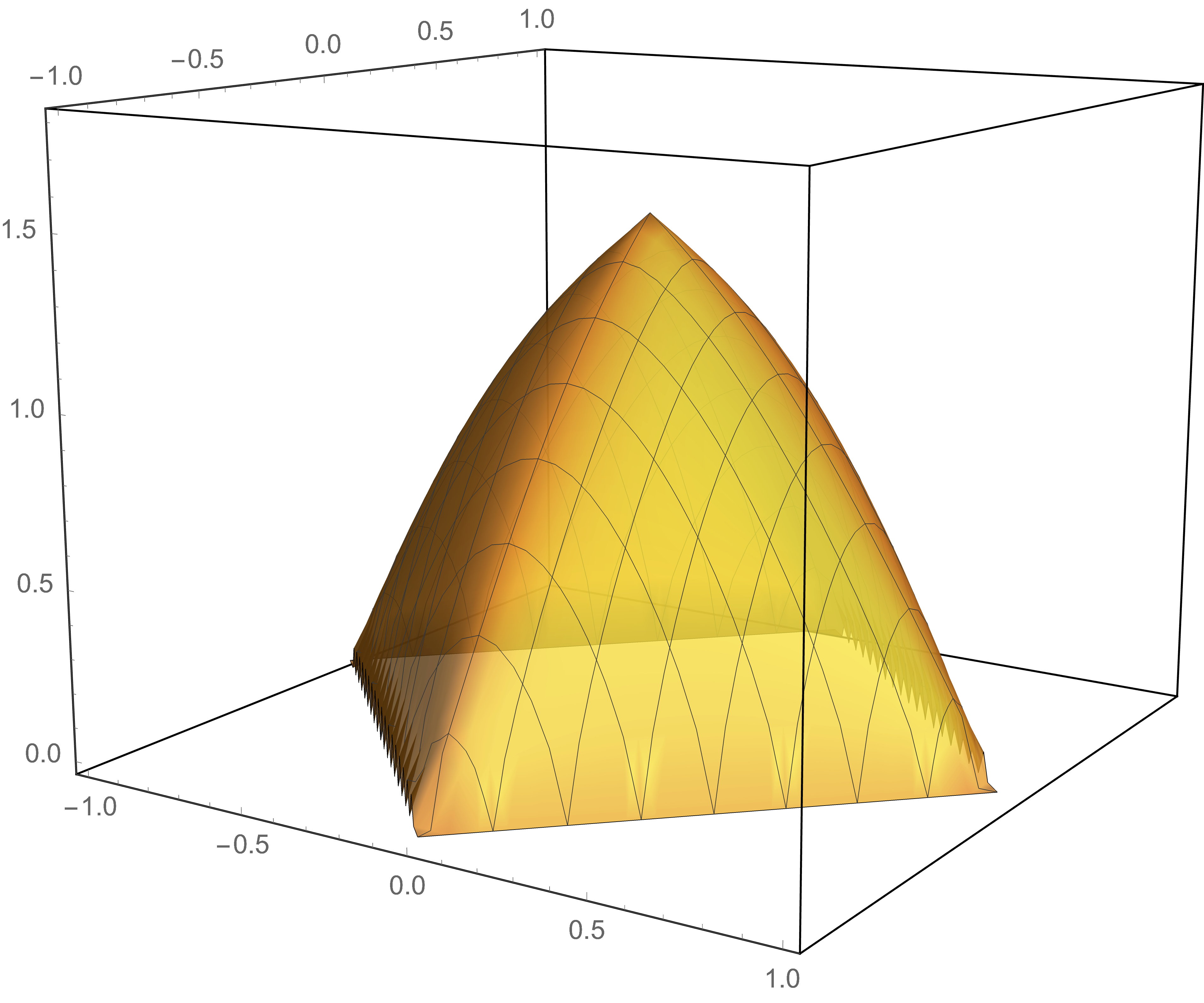}\caption{}
   	\end{subfigure}
   	\hfil
		\begin{subfigure}[b]{0.38\textwidth}
			\centering
	   	\includegraphics[width=\textwidth]{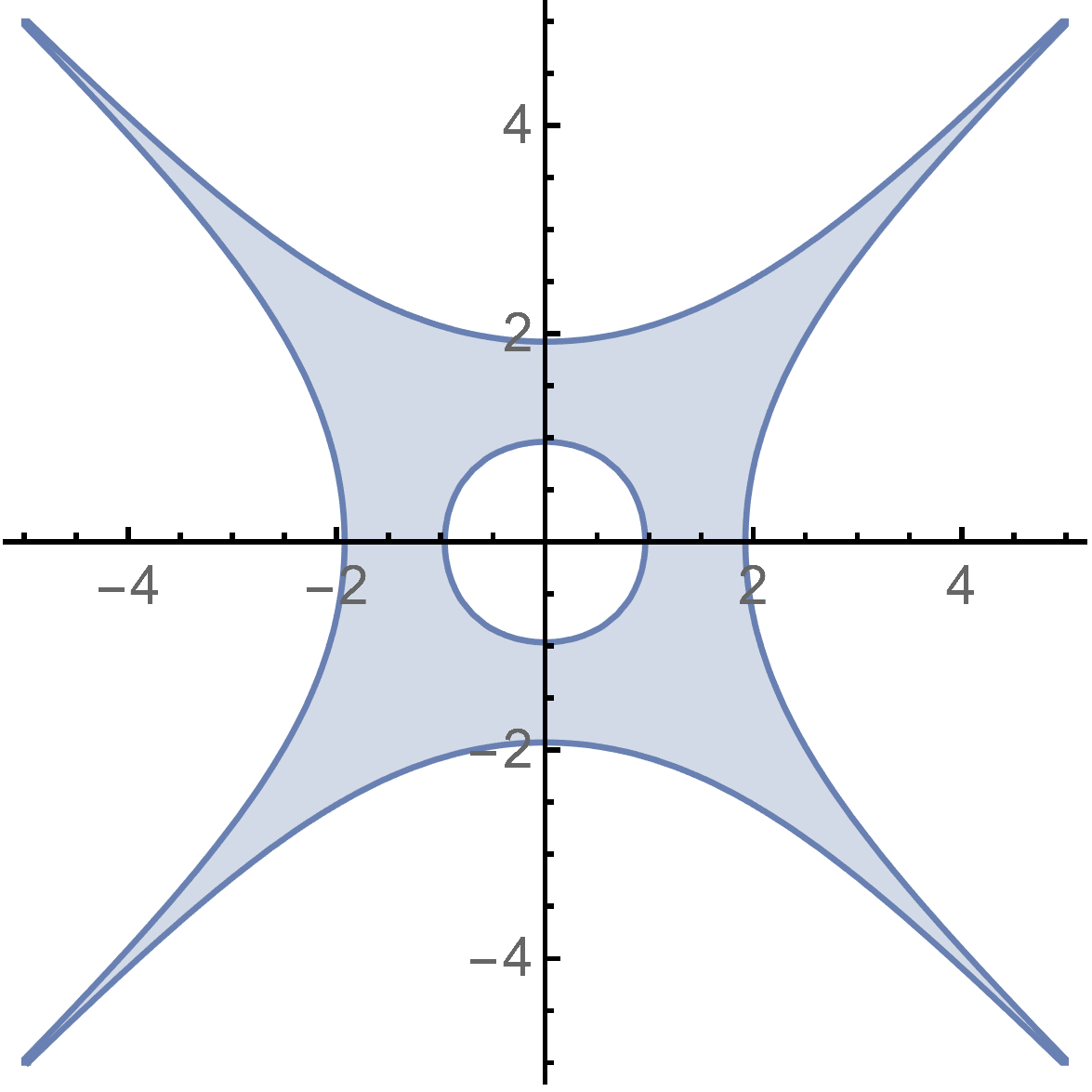}\caption{}
		\end{subfigure}
		\caption{(a) The decimation limit for $f(x,y)= 5 + x + x^{-1}+y+y^{-1}$ from Example
		      \ref{exam:amoeba-with-hole} \label{fig:amoeba-with-hole}, and (b) the ``hole'' in its amoeba causing the
		      peak.}
	\end{figure}

   As in the previous example, the decimation limit describes the surface tension for a
   physical model, in this case dimer tilings of the square-octagon graph (see
  \cite{KenyonOkounkovSheffield}*{Fig.\ ~3}. 
\end{example}

\begin{remark}\label{rem:1+x+y}
	Dimer models have a long history in statistical physics. A particularly important instance involves $f(x,y)=1+x+y$ from Example \ref{exam:ledrappier}, and has been studied in enormous detail by many authors, including Kenyon, Okounkov, and Sheffield \cite{KenyonOkounkovSheffield}.
	
	To describe this model, let $\mathscr{H}$ denote the regular hexagonal lattice in $\RR^2$. We can assign the vertices of $\mathscr{H}$ alternating colors red and black, much like a checkerboard. A \emph{perfect matching} on $\mathscr{H}$ is an assignment of each red vertex to a unique adjacent black vertex, these forming an edge or \emph{dimer}. A perfect matching is equivalent to a tiling of $\RR^2$ by three types of lozenges, one type for each of the three edges incident to each vertex. Using a natural height function, such a lozenge tiling gives a surface, and the study of the statistical properties of such random surfaces has resulted in many remarkable discoveries (see Okounkov's survey \cite{Okounkov} or Gorin's detailed account of lozenge tilings \cite{Gorin}).
	
	 discovered that by cleverly assigning signs to the edges of $\mathscr{H}$, he could compute the number of perfect matchings on a finite approximation using periodic boundary conditions by a determinant formula. Furthermore, this determinant can be explicitly evaluated to have the form of a decimation of $f(x,y)=1+x+y$. Each of the three terms of $f$ correspond to one of the three types of lozenges in the random tiling. Then according to \cite{KenyonOkounkovSheffield}*{\S3.2} the logarithmic scaling limit $\df(r,s)$ counts the growth rate of perfect matchings for which the frequencies of the three lozenge types are $r$, $s$, and $1-r-s$. As such, it is called the \emph{surface tension} for this model.
	
	The two-variable polynomials with integer coefficients arising from such dimer models, such as the preceding two examples,  define curves of a very special type called Harnack curves. For these there are probabilistic interpretations of the coefficients of decimations. The additional structure enables one to  show that the individual nonzero coefficients of $\fn$ grow at a rate predicted by $\df$.  Example \ref{exam:transcendental} shows this can fail if complex coefficients are allowed. But whether or not this is true for every polynomial in $\zzd$ for all $d\ge1$ appears to be quite an interesting problem (see Question \ref{que:coefficcients} for a precise formulation).
\end{remark}

\section{Convex functions and Legendre duals}\label{sec:legendre-duals}

We briefly review some basic facts about convex functions and their Legendre duals. Rockafellar's classic book \cite{Rockafellar} contains a comprehensive account of this theory.

Let $\oR$ denote $\RR\cup\{\infty\}$, with the standard conventions about arithmetic operations
and inequalities involving $\infty$. Let $\phi\colon \rd\to\oR$ be a function, and define its \emph{epigraph} by
\begin{displaymath}
   \epi\phi :=\{ (\bu,t): \bu\in\rd, t\in\RR, \text{\ and\ } t\ge \phi(\bu) \}\subset \rd\times\RR.
\end{displaymath}
Then $\phi$ is  defined to be \emph{convex} provided that $\epi\phi$ is a convex subset of $\rd\times \RR$. A~ function $\psi\colon\rd \to\uR$ is called \emph{concave} if $-\psi\colon\rd\to\oR$
is convex.

The \emph{effective domain} of a convex function $\phi$ is defined by
\begin{displaymath}
   \dom\phi :=\{\bu\in\rd:\phi(\bu)<\infty\}. 
\end{displaymath}
By allowing $\phi$ to take the value $\infty$, we may assume that it is defined on all of
$\rd$, enabling us to combine convex functions without needing to take into account their
effective domains. A convex function is \emph{closed} if its epigraph is a closed subset of
$\rd\times\RR$. This property normalizes the behavior of a convex function at the boundary of
its effective domain, and holds for all convex (and concave) functions that arise here.

Suppose that $\phi\colon\RR\to\oR$ is convex. Its \emph{Legendre dual} (or, more accurately, its \emph{Legendre--Fenchel dual}) $\phi^*$ is defined for all $\br\in\rd$ by
\begin{equation}\label{eqn:dual}
   \phi^*(\br):= \sup\{ \<\br,\bu\>-\phi(\bu) : \bu\in\rd\}.
\end{equation}
The Legendre dual $\phi^*$ is also a convex function, and provides an alternative description of $\epi\phi$ in terms of its support hyperplanes. Furthermore, Legendre duality states that $\phi^{**}=\phi$ for closed convex functions.

The Legendre dual of a concave function $\psi\colon\rd\to\uR$ is similarly defined as
\begin{equation}\label{eqn:concave-dual}
   \psi^*(\br)=\inf\{\<\br,\bu\>-\psi(\bu):\bu \in\rd\}.
\end{equation}
Then $\phi=-\psi$ is convex, and a simple manipulation shows that their Legendre duals are related by $\psi^*(\br)=-\phi^*(-\br)$.

\section{Amoebas and Ronkin functions}\label{sec:ronkin-functions}

Let $0\ne f\in\czd$. Put $\CC^*=\CC\smallsetminus\{0\}$ and define $\variety(f):=
\{\bz\in(\CC^*)^d:f(\bz)=0\}$. Let $\operatorname{Log}\colon(\CC^*)^d\to\rd$ be the map
$\operatorname{Log}(z_1,\dots,z_d)=(\log|z_1|,\dots,\log |z_d|)$.

\begin{figure}[t!]
    \centering
    \begin{subfigure}[b]{0.45\textwidth}
    	\centering
    	\includegraphics[width=0.8\textwidth]{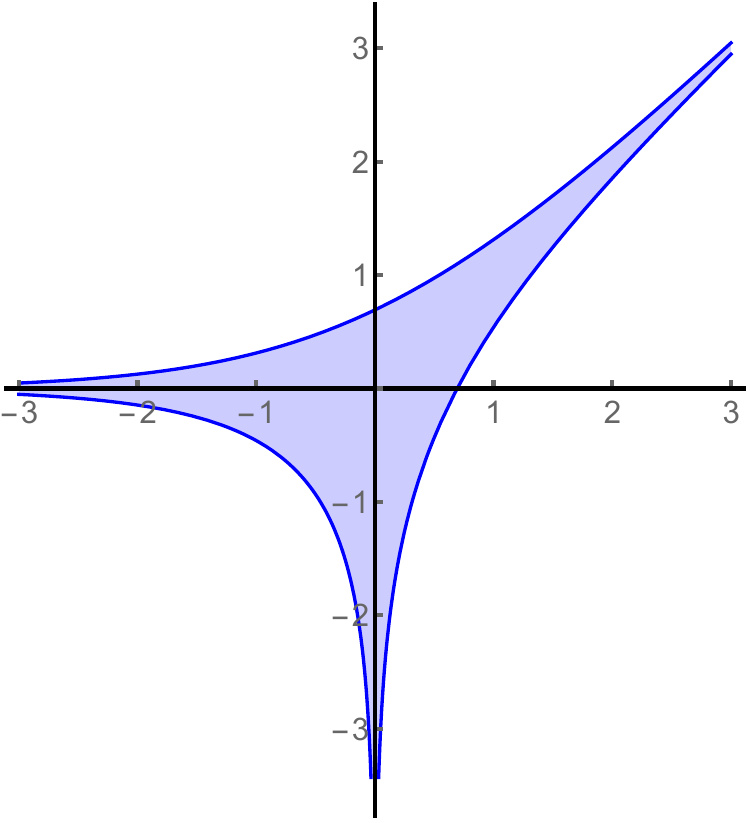}
    	\caption{}
    \end{subfigure}
    \hfil
    \begin{subfigure}[b]{0.45\textwidth}
      \centering
      \includegraphics[width=\textwidth]{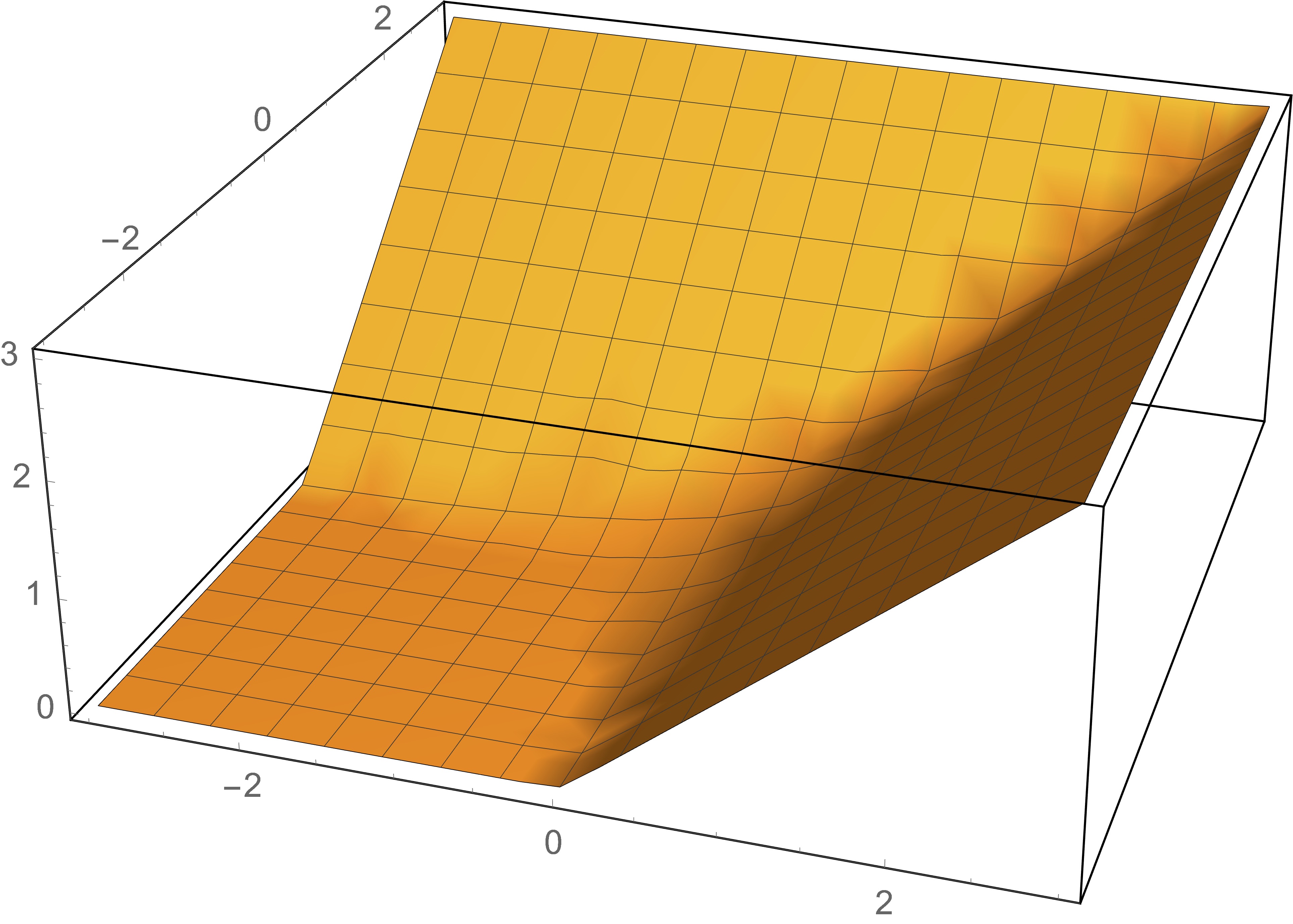}
      \caption{}
    \end{subfigure}
    \caption{(a) The amoeba of $1+x+y$, and (b) its Ronkin function \label{fig:amoeba}}
\end{figure}

In 1993 Gelfand, Kapranov, and Zelevinsky \cite{GKZ} introduced the notion of the \emph{amoeba}
$\SA_f$ of $f$, defined as 
\begin{displaymath}
   \SA_f:= \operatorname{Log}\bigl( \variety(f) \bigr) \subset \rd.
\end{displaymath}
The amoeba of $1+x+y$ is depicted in Figure \ref{fig:amoeba}(a). The complement
$\SA_f^c=\rd\smallsetminus \SA_f^{}$ of $\SA_f$ consists of a finite number of connected
components, all convex. The unbounded components are created by ``tentacles'' of $\SA_f$.
Unfortunately, biological amoebas look nothing like their mathematical namesakes.

Closely related to $\SA_f$ is the Ronkin function $\rf$ of $f$, introduced by Ronkin
\cite{Ronkin} in 2001, and defined earlier in \eqref{eqn:ronkin}. The Ronkin function of $1+x+y$ is shown in Figure \ref{fig:amoeba}(b).

The Ronkin function of a polynomial $f$ is known to be a convex function on $\rd$ and affine on
each connected component of $\SA_f^c$ (see
\cite{PassareRullgard} for all properties of $\ssR_f$ and $\SA_f$ used here). Moreover, on each connected component of $\SA_f^c$ the (constant)
gradient of $\rf$ is contained in $\nf\cap\zd$, and the convex hull of these values equals
~$\nf$. From this we conclude that the Legendre dual $\ssR_f^*$ of $\ssR_f^{}$ has effective
domain ~$\nf$.

\section{Decimation limits of polynomials}\label{sec:full-decimation}

In this section we prove Theorem \ref{thm:main}, one of our main results, and Corollaries~ \ref{cor:max} and~ \ref{cor:tropical-convolution}.  If $0\ne f\in\czd$ we will show that the $N$th renormalized decimation $\dnf=\ch(\LN f)$ converges uniformly on $\nf$ to a continuous concave limit function $\df$, and that $\dlf=-\rf^*$.

The first ingredient in our proof is the basic estimate of Mahler relating the
largest coefficient of a polynomial to its Mahler measure and its support. Let us begin with some terminology. For $0\ne g \in \czd$ define its \emph{height} $\Ht(g)$ by
$\Ht(g)=\max\{\,|\gh(\bk)|:\bk\in\zd\}$. The \emph{Mahler measure} of $g$ is
$\Mahler(g)=\exp\bigl(\mahler(g)\bigr)$, where $\mahler(g)$ is the logarithmic Mahler measure defined in~
\eqref{eqn:mahler-measure}. 

\begin{proposition}[Mahler \cite{Mahler}]\label{prop:mahler}
   Suppose that $0\ne g\in\czd$ and that $\supp g\subset [0,C-1]^d\cap\zd$. Then
   \begin{equation}\label{eqn:mahler-inequality}
      2^{-dC} \Ht(g) \le \Mahler(g) \le C^d \Ht(g).
   \end{equation}
\end{proposition}

\begin{proof}
   Let $\bk=(k_1,\dots,k_d)\in\supp g$. Then by \cite{Mahler}*{Eqn.\ (3)},
   \begin{displaymath}
      |\gh(\bk)|\le \binom{C-1}{k_1} \binom{C-1}{k_2}\dots \binom{C-1}{k_d} \Mahler(g).
   \end{displaymath}
   Since each binomial coefficient is bounded above by $2^C$, the first inequality in
   \eqref{eqn:mahler-inequality} follows.

   To prove the second inequality, observe that for all real numbers $s_1,\dots,s_d$ we have that
   \begin{equation}\label{eqn:length-inequality}
      |g(e^{2\pi i s_1},\dots,e^{2\pi i s_d})| \le \sum_{\bk\in\zd} |\gh(\bk)| \le
      |[0,C-1]^d\cap\zd|\cdot\Ht(g)
      = C^d\Ht(g). 
   \end{equation}
   Hence
   \begin{displaymath}
      \Mahler(g) = \exp\Bigl[ \int_0^1 \dots \int_0^1 \log|g(e^{2\pi i s_1},\dots, e^{2\pi i s_d})| \,ds_1\dots ds_d \Bigr] \le C^d\Ht(g). \qedhere
   \end{displaymath}   
\end{proof}

\qedhere

Consider $\cstd$ as a group under coordinate-wise multiplication. Define the action of $\bz\in\cstd$ on $f\in\czd$ by $(\bz\cdot f)(x_1,\dots,x_d)=f(z_1x_1,\dots,z_dx_d)$. This action is commutative since
\begin{displaymath}
   \bz\cdot(\bz'\cdot f) = (\bz \bz')\cdot f = \bz'\cdot(\bz\cdot f),
\end{displaymath}
and also $\bz\cdot(fg)=(\bz\cdot f)(\bz\cdot g)$ for all $f,g\in\czd$. Hence the map $f\mapsto \bz\cdot f$ is a ring isomorphism of $\czd$. Furthermore, $(\bz\cdot f)\sphat\,(\bk)=\bz^{\bk} \fh(\bk)$ for all $\bk\in\zd$, and so $\SN_{\bz\cdot
f}=\nf$ for all $\bz\in\cstd$.

Recall that $\Omn$ denotes the group of $N$th roots of unity. For $\bo\in\Omnd\subset\cstd$ we call $\bo\cdot f$ the \emph{rotate of $f$ by $\bo$}. Then $\fn=\prod_{ \bo\in\Omnd} \bo\cdot f$ is the product of all rotates of $f$ by elements in $\Omnd$.

If $g,h\in\czd$ then it is well known that $\SN_{gh}=\SN_g+\SN_h$ (the Minkowski sum), and
trivially $\ssR_{gh}=\ssR_g+\ssR_h$. By our previous remarks,
\begin{displaymath}
   \SN_{\fn}=\sum_{\bo\in\Omnd}\SN_{\bo\cdot f}=\sum_{\bo\in\Omnd}\nf= N^d\nf.
\end{displaymath}
Also, $\ssR_{\bo\cdot f}=\rf$, and hence $\ssR_{\fn}=N^d\ssR_{f_{}}$.
 
For $\bu\in\rd$ put $\eu=(e^{u_1},\dots,e^{u_d})$. Then $(\eu\cdot f)\widehat{\
}(\bk)=e^{\bu\cdot\bk}\fh(\bk)$. Commutativity of the action of $\cstd$ on $f$ then shows that
$(\eu\cdot f)_{\<N\>}=\eu\cdot(\fn)$. Also
\begin{displaymath}
   \rf(\bu)=\log\Mahler(\eu\cdot f)=\frac{1}{N^d}\log\Mahler\bigl( (\eu\cdot f)_{\<N\>}\bigr) =
   \frac{1}{N^d} \log \Mahler(\eu\cdot\fn).
\end{displaymath}
Observe that 
\begin{displaymath}
   \log \Ht(\eu\cdot \fn)=\max\{\<\bu,\bk\>+\log|\fnh(\bk)|:\bk\in\zd\},
\end{displaymath}
indicating a connection with Legendre duals.

\begin{proof}[Proof of Theorem \ref{thm:main}]
   Let $0\ne f\in\czd$. Fix $\bfm\in\zd$ and let $g(\bx)=\bx^{\bfm}f(\bx)$. It is straightforward to
   verify that $(\DN g)(\br)=(\DN f)(\br-\bfm)$ for all $\br\in\rd$. Therefore by adjusting $f$
   by suitable monomial, we may assume that $\supp f\subset[0,B-1]^d\cap\zd$ for some $B\ge
   1$. Then $\supp (\eu\cdot\fn)\subset [0,N^d(B-1)]^d\cap\zd\subset[0,N^dB-1]^d\cap\zd$ for
   every $\bu\in\rd$. By Proposition \ref{prop:mahler},
   \begin{align*}
     \rf(\bu) &= \frac{1}{N^d} \log\Mahler(\eu\cdot\fn) \le \frac{1}{N^d} \Bigl\{
                \log\bigl[(N^dB)^d\bigr] +\log \Ht(\eu\cdot\fn) \Bigr\} \\
               &= \frac{\log\bigl[(N^dB)^d\bigr]}{N^d}+ \frac{1}{N^d}\max_{\bk\in\zd}
     \{\<\bu,\bk\> + \log |\fnh(\bk)|\},
   \end{align*}
   where the error term $b_N:=N^{-d}\log\bigl[(N^dB)^d\bigr]\to 0$ as $N\to\infty$, uniformly
   for $\bu\in\rd$.

   An opposite inequality is based on the following fundamental observation, used both by Boyd
   \cite{BoydUniform} and Purbhoo \cite{Purbhoo} for different purposes. As we noticed before, $\fn$
   is a polynomial in the $N$th powers of the variables. Therefore $\EN\fnh$ is again a
   polynomial to which we can apply Prop.\ \ref{prop:mahler}, but with improved constants since
   the support has now shrunk by a factor of $N$. This improvement is crucial.

   Specifically,
   \begin{displaymath}
      \supp (\eu\cdot\fn)\subset [0,N^d(B-1)]^d\cap (N\zd),
   \end{displaymath}
   so that
   \begin{displaymath}
      \supp(\EN(\eu\cdot\fn))\subset[0,N^{d-1}(B-1)]\cap\zd.
   \end{displaymath}
   Applying Prop.\ \ref{prop:mahler},
   \begin{equation}\label{eqn:Mahler-bound}
   	\begin{aligned}
   	   \Ht(\eu\cdot\fn)&=\Ht(\EN(\eu\cdot\fn))\le 2^{dN^{d-1}B}\Mahler(\EN(\eu\cdot\fn))\\
   	   &=2^{dN^{d-1}B}\Mahler (\eu\cdot f)^{N^d}.
   	\end{aligned}
   \end{equation}
   Hence
   \begin{displaymath}
      \frac{1}{N^d} \log \Ht(\eu\cdot\fn)\le \frac{dN^{d-1}B\log 2}{N^d} +\log\Mahler(\eu\cdot
      f) = a_N+\rf(\bu),
   \end{displaymath}
   where again the error term $a_N:= (d B\log 2)/N\to0$ uniformly for $\bu\in\rd$.
   We can summarize these estimates as
   \begin{equation}\label{eqn:uniform-bounds}
      \Bigl|\rf(\bu) - \frac{1}{N^d} \max_{\bk\in\zd}\bigl\{\<\bu,\bk\>+\log|\fnh(\bk)|\bigr\}\Bigr|\le
      \max\{a_N,b_N\}\to 0
   \end{equation}
   as $N\to\infty$ uniformly in $\bu\in\rd$.

   Next we relate the first max occurring in \eqref{eqn:uniform-bounds} with the $N$th normalized decimation $\dnf$. We have that
   \begin{align*}
     \frac{1}{N^d} \max_{\bk\in\zd} &  \{\<\bu,\bk\>+\log|\fnh(\bk)|\} =
          \max_{\bk\in\zd} \Bigl\{ \Bigl<\bu,\Bigl(\frac{\bk}{N^d}\Bigr)\Bigr> +
                                      \frac{1}{N^d}\log|\fnh(\bk)| \Bigr\}\\
     & = \max_{\bk\in\zd} \Bigl\{ \Bigl<\bu,\Bigl(\frac{\bk}{N^d}\Bigr)\Bigr> +
       \frac{1}{N^d}\ssE_{N^d}\log \Bigl|\fnh\Bigl(\frac{1}{N^d}\bk\Bigr)\Bigr|\Bigr\}\\
       &= \max_{\bk\in\zd} \Bigl\{ \Bigl<\bu,\Bigl(\frac{\bk}{N^d}\Bigr)\Bigr> +
       (\dnf)\Bigl(\frac{\bk}{N^d}\Bigr)\Bigr\}\\
       & = \max_{\br\in\rd} \{\<\bu,\br\> + \dnf(\br)\} = -(\dnf)^*(-\bu).
   \end{align*}

   Hence by \eqref{eqn:uniform-bounds}, $-(\dnf)^*(-\bu)$ converges to $\rf(\bu)$ uniformly for
   $\bu\in\rd$, or, equivalently,
   \begin{equation}\label{eqn:convergence}
     (\dnf)^*(\bu)\to -\rf(-\bu) \text{\quad uniformly for $\bu\in\rd$}.
   \end{equation}

   If $\phi$ and $\psi$ are concave functions on $\rd$ such that
   $|\phi(\bu)-\psi(\bu)|\le \varepsilon$ for all $\bu\in\rd$, it is easy to check from the
   definitions that $\phi^*$ and $\psi^*$ have the same effective domain, and that
   $|\phi^*(\br)-\psi^*(\br)|\le\varepsilon $ for all $\br\in\dom \phi^*=\dom\psi^*$. Applying
   this to \eqref{eqn:convergence} and using duality we finally obtain that
   $(\dnf)^{**}=\dnf\to -\rf^*$ uniformly on $\nf$, completing the proof.
\end{proof}

\begin{proof}[Proof of Cor.\ \ref{cor:max}:] 
   By Theorem \ref{thm:main}, Legendre duality, and \eqref{eqn:concave-dual},
   \begin{displaymath}
        -\mahler(f)=-\rf(0,0)=\ssD_f^{*}(0,0) = \inf_{(r,s)\in\nf} -\df(r,s) = -\mspace{-6mu} \sup_{(r,s)\in\nf} \df(r,s). 
   \end{displaymath}
   We remark that differentiability of $\df$ at the maximum value is not assumed for Legendre duality to apply here, and Example \ref{exam:amoeba-with-hole} provides a case when differentiability fails.
\end{proof}

\begin{proof}[Proof of Cor.\  \ref{cor:tropical-convolution}:]
   Let $f$ and $g$ be nonzero polynomials in $\czd$. Clearly $\ssR_{fg}=\ssR_f + \ssR_g$. By
   \cite{Rockafellar}*{Thm.\ 16.4}, the
   Legendre dual of the sum $\phi+\psi$ of two convex functions is their infimal convolution
   defined for $\br\in\rd$ by $\inf\{\phi(\bs) + \psi(\br-\bs):\bs\in\rd\}$. Applying
   this with $\phi = -\ssR_f$ and $\psi=-\ssR_g$, using Thm.\ \ref{thm:main}, and taking
   negatives  we obtain that $\ssD_{fg}=\ssD_f\circledast\ssD_g$.
\end{proof}

\begin{remark}
   Our estimate \eqref{eqn:uniform-bounds} can be expressed in the language of tropicalization of polynomials (see \cite{MS}*{\S3.1} for background and motivation). Let $0\neq g(\bx)= \sum_{\bk\in\zd} \gh(\bk)\bx^{\bk}\in\czd$. Define the \emph{tropicalization} of $g$ to be the function $\trop g\colon \RR^d\to\RR$ given by
   \begin{displaymath}
      (\trop g)(\bu) = \max_{\bk\in\zd} \bigl\{ \<\bu,\bk\>+\log|\gh(\bk)|\bigr\},
   \end{displaymath}
   which is a polyhedral convex function.
   Then by \eqref{eqn:uniform-bounds} we see that
   \begin{equation}\label{eqn:tropicalization}
      \frac{1}{N^d} \trop \fn \to \rf \text{\quad uniformly on $\RR^d$,}
   \end{equation}
   so that the normalized tropicalization of $\fn$ converges uniformly to the Ronkin function of~$f$. Figure \ref{fig:polyhedralronkin}(a) depicts this polyhedral approximation for $f(x,y)=1+x+y$ and $N=5$ (compare with Figure \ref{fig:amoeba}(b)). The tropical variety of this polyhedral approximation is the projection to the plane of the vertices and edges of its graph, and is shown in \ref{fig:polyhedralronkin}(b). These tropical varieties converge in the Hausdorff metric to the amoeba of ~$f$ as $N\to\infty$ (compare with Figure \ref{fig:amoeba}(a)).

   \begin{figure}[ht]
       \centering
       \begin{subfigure}[b]{0.45\textwidth}
       	\centering
       	\includegraphics[width=1.1\textwidth]{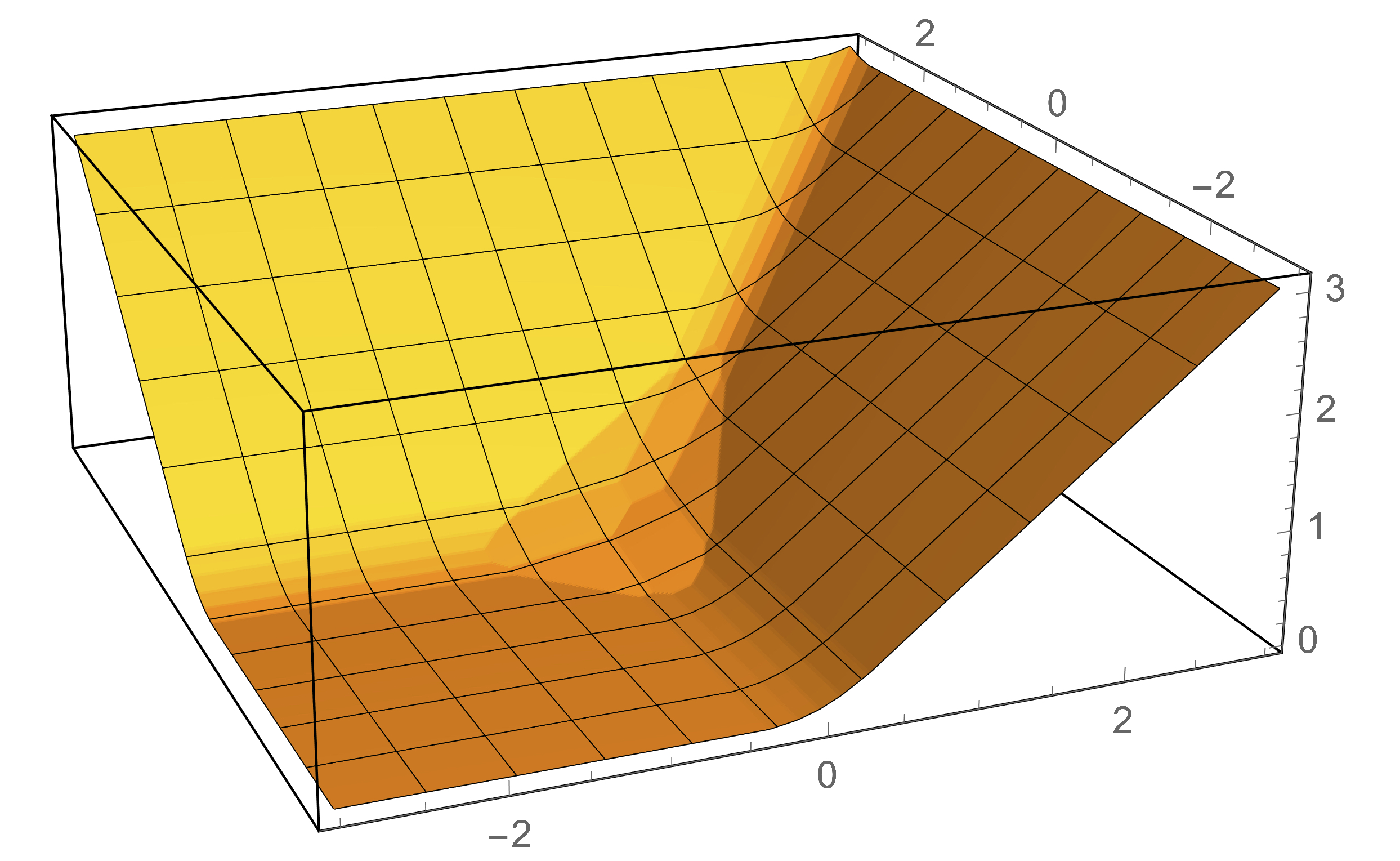}
       	\caption{}
       \end{subfigure}
       \hfil
       \begin{subfigure}[b]{0.35\textwidth}
         \centering
         \includegraphics[width=\textwidth]{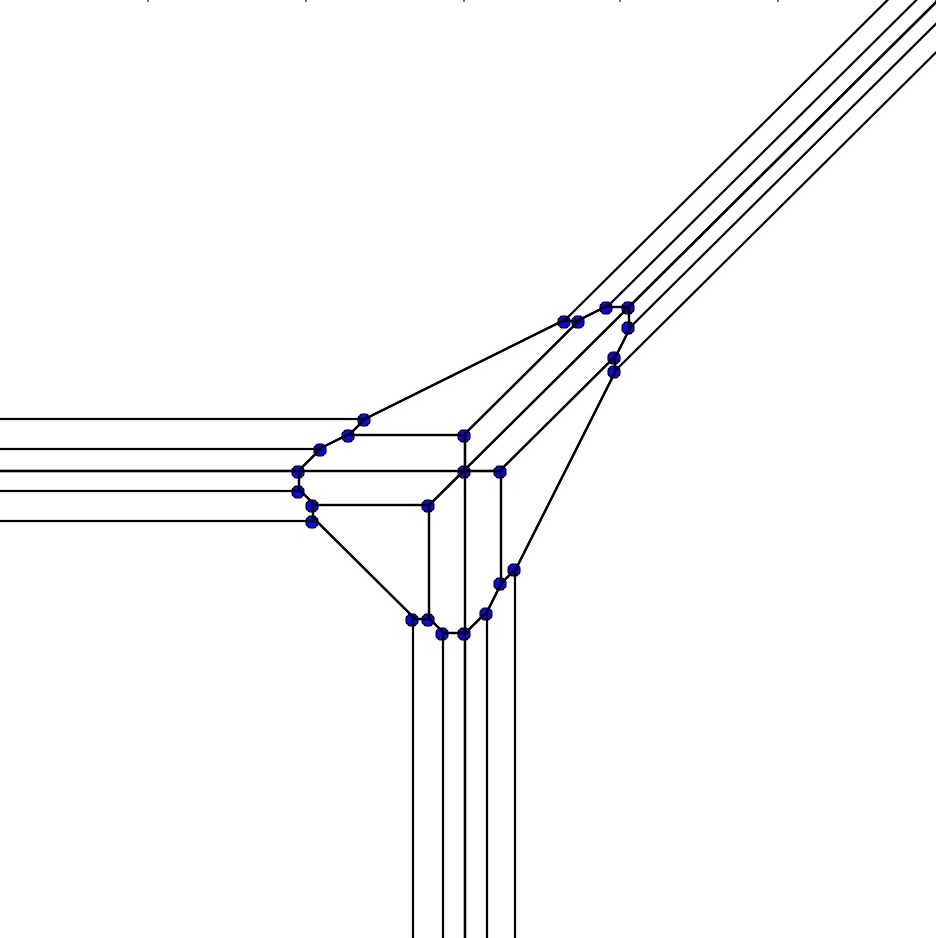}
         \caption{}
       \end{subfigure}
       \caption{(a) Tropical approximation to the Ronkin function of $1+x+y$, and (b) its corresponding tropical variety \label{fig:polyhedralronkin}}
   \end{figure}   
\end{remark}

\begin{remark}
	In \cite{Purbhoo} Purbhoo used decimations for a different purpose, namely to find a computational way to detect whether or not a point is in the amoeba of a given polynomial. 
	Call a polynomial \emph{lopsided} if it has one coefficient whose absolute value strictly exceeds the sum of the absolute values of all the other coefficients. Let $f\in\zzd$ and $\bu\in\rd$. Clearly if $e^{\bu}\cdot f$ is lopsided then $\bu\notin\Af$. Purbhoo used decimations to amplify size differences among the coefficients. More precisely,
	he proves that given $\varepsilon>0$ there is an $N_0$, depending only on $\varepsilon$ and the support of $f$, such that if $N>N_0$ and the distance from $\bu$ to $\Af$ is greater than $\varepsilon$ then $e^{\bu}\cdot\fn$ is lopsided. Since $f$ and $\fn$ have the same amoeba, this gives an effective algorithm for approximating the complement of $\Af$.
	
	One direct consequence of \cite{Purbhoo} is that the normalized tropicalizations in \eqref{eqn:tropicalization} converge to the Ronkin function off the amoeba of $f$, while our result is that this convergence is uniform on all of $\rd$. Roughly speaking, Purbhoo is concerned with the coefficients of $e^{\bu}\cdot f$ for points $\bu$ off the amoeba, while our focus is on $\bu$ within the amoeba.
\end{remark}

\begin{remark}
   Let $F$ be a lower-dimensional face of the Newton polytope $\nf$ of ~$f$, and put $f|_F=\sum_{\bn\in F} \fh(\bn)\bx^{\bn}$. Clearly the restriction of $\df$ to $F$ is just the decimation limit of $f|_F$, or in symbols $\df|_F=\ssD_{f|_F}$. By Corollary \ref{cor:max}, this generalizes \cite{LSW}*{Rem.\ 5.5}, which gave a dynamical proof of the inequality due to Smyth \cite{SmythKronecker}*{Thm.\ 2} that $\mahler(f)\ge\mahler(f_F)$ for every face $F$ of $\nf$.
\end{remark}

\section{Decimations of principal actions and contracted
ideals}\label{sec:algebraic-decimations}

We return to decimations of principal algebraic $\zd$-actions, and in this section show that
they are again principal. The proof uses machinery from commutative algebra, including
contractions of ideals.

Suppose that $X$ is a compact, shift-invariant subgroup of $\tzd$. Using Pontryagin duality we
can obtain an alternative description of $X$ as follows (for a comprehensive account see
\cite{SchmidtBook}*{Chap.\ II}).

As a discrete abelian group the Pontryagin dual of $\tzd$ is the direct sum of $\zd$ copies of
$\ZZ$, which we suggestively write as $\bigoplus_{\bk\in\zd}\ZZ \bx^{\bk}=\zzd$. The (additive)
dual pairing between $\tzd$ and $\zzd$ is given by $\< t,g \>=\sum_{\bk\in\zd}
t_{\bk}\,\gh(\bk)\in\TT$. Multiplication by the inverses of each of the variables $x_j$ on $\zzd$ gives a $\zd$-action that is dual to the natural shift action $\sigma$ on $\tzd$ defined earlier. 

Since $X$ is shift-invariant, $\{g\in\zzd : \<t,g\>=0 \text{\  for all $t\in X$}\}$ is an
ideal $\fa$ in $\zzd$, and the dual group of $X$ equals $\zzd/\fa$. Conversely, if $\fa$ is an
arbitrary ideal in $\zzd$, then the compact dual group $X_{\fa}$ of $\zzd/\fa$ is a
shift-invariant subgroup of $\tzd$. Thus there is a one-to-one correspondence between
shift-invariant compact subgroups of $\tzd$ and ideals in $\zzd$. When $\fa$ is the principal
ideal $\id{f}$ generated by ~$f$, then $X_{\fa}=\xf$ as defined above, explaining the
terminology ``principal actions''.

Fix $N\ge 1$ and recall the restriction map $r_N\colon \tzd\to\tnzd$ from
\S\ref{sec:introduction}. Let $f\in\zzd$. Then the $N$th decimation $r_N(\xf)$ is a compact
subgroup of $\tnzd$ that is invariant under the shift-action of $\nzd$. By our previous
discussion, the dual group of $r_N(\xf)$ has the form $\znzd/\fa_N$, where $\fa_N$ is an ideal
in $\znzd$. The following result identifies this ideal.

\begin{lemma}\label{lem:contracted-ideal}
   Let $f\in\zzd$ and $N\ge 1$. Then the dual group of $r_N(\xf)$ is $\znzd/\fa_N$, where
   $\fa_N=\id{f}\cap\znzd$. 
\end{lemma}

\begin{proof}
   Let $\fb_N=\{g\in\znzd:\<t,g\> = 0 \text{\ for all $t\in r_N(\xf)$}\}$. If $g\in\fa_N$, then
   for every $t\in\xf$ we have that $0=\< t,g\>=\<r_N(t),g\>$, so that $g\in\fb_N$. Conversely,
   if $g\in\fb_N$ and $t\in \xf$, then $g$ annihilates the restriction of $t$ to every coset of
   $\nzd$, and hence annihilates $t$, so that $g\in\fa_N$.
\end{proof}

The ideal $\id{f}\cap\znzd$ defining $r_N(\xf)$ is called the \emph{contraction} of $\id{f}$ to $\znzd$. The main result of this section is that this contraction is always principal.

\begin{proposition}\label{prop:contractions-are-principal}
   Let $f\in\zzd$ and $N\ge 1$. Then the contracted ideal $\id{f}\cap\znzd$ is a principal
   ideal in $\znzd$. 
\end{proposition}

We begin by briefly sketching the necessary terminology and machinery from commutative algebra,
all of which is contained in \cite{AtiyahMacdonald} or can be easily deduced from material there.

For brevity let $R=\znzd$ and $S=\zzd$. Both $R$ and $S$ are unique factorization domains, and
therefore both are integrally closed \cite{AtiyahMacdonald}*{Prop.\ 5.12}. Furthermore, $S$ is integral over $R$ since each variable $x_j$ in $S$ satisfies the monic polynomial $y_{}^N-x_j^N\in R[y]$.

A prime ideal $\fp$ in an integral domain has \emph{height one} if there are no prime ideals
strictly between $0$ and $\fp$. In a unique factorization domain the prime ideals of height one
are exactly the principal ideals generated by irreducible elements. A proper ideal ~$\fq$ in an
integral domain is \emph{primary} if whenever $ab\in\fq$ then either $a\in\fq$ or $b^n\in\fq$
for some $n\ge1$. In this case its radical $\{a:a^n\in\fq\text{\ for some $n\ge1$}\}$ is a
prime ideal, say~ $\fp$, and then $\fq$ is called \emph{$\fp$-primary}. Examples show that in general a power of a prime ideal need not be primary, that a primary ideal need not be the power of a prime ideal, and that even if an ideal has prime radical it need not be primary. The notion of primary ideal, although the correct one for decomposition theory, is quite subtle. However, in our situation things are much simpler.

\begin{lemma}\label{lem:primary}
   Let $P$ be a unique factorization domain, and let $r\in P$ be irreducible. Then the principal ideal $\fp=\id{r}$ is prime, and the $\fp$-primary ideals are exactly the powers $\fp^n$ of $\fp$ for $n\ge 1$.
\end{lemma}

\begin{proof}
   It is clear that $\fp$ is prime. To prove that $\fp^n=\id{r^n}$ is $\fp$-primary, suppose that $ab\in\fp^n$, but $a\notin \fp^n$. Then $r\mid b$, so $b^n\in\fp^n$, showing that $\fp^n$ is primary. Clearly the radical of $\fp^n$ is $\fp$, and so $\fp^n$ is $\fp$-primary.
   
   Conversely, suppose that $\fq$ is a $\fp$-primary ideal. Since the radical of $\fq$ is $\fp$, it follows that $r^n\in\fq$ for some $n\ge1$. Choose $n$ to be the minimal such power, so that $\fp^n\subset\fq$. Suppose that $\fp^n \neq \fq$, and let $a\in \fq\smallsetminus\fp^{n}$. Write $a=cr^m$, where $r\nmid c$. Clearly $m\le n-1$, and so $r^m\notin\fq$ by minimality of $n$. Since $\fq$ is primary, there is a $k\ge1$ such that $c^k\in\fq\subset\fp$. But this contradicts $r\nmid c$. Hence $\fp^n = \fq$.
\end{proof}

If $\fa$ is an ideal in $S$, we denote its \emph{contraction} $\fa\cap R$ to $R$ by
$\fa\sfc$. If $\fq$ is a $\fp$-primary ideal in $S$, then $\fp\sfc$ is prime and
$\fq\sfc$ is $\fp\sfc$-primary in $R$.

One of the important results in commutative algebra, essential to developing a dimension theory
using chains of prime ideals, is the so-called ``Going Down'' theorem
\cite{AtiyahMacdonald}*{Thm.\ 5.16}.  Its hypotheses are satisfied in our situation, and it
says the following. Suppose that $\fp_0\subsetneq \fp_1 \subsetneq \fp_2$ is a chain of prime
ideals in $R$, and that there is a prime ideal $\fq_2$ in $S$ with
$\fq_2\sfc=\fp_2^{}$. Then there is a chain $\fq_0\subsetneq \fq_1 \subsetneq \fq_2$ of
prime ideals in $S$ such that $\fq_j\sfc=\fp_j^{}$ for $j=0,1,2$. From this it follows that
prime ideals in $S$ of height one contract to prime ideals in $R$ of height one. In other
words, if $h\in S$ is irreducible, then $\ids{h}\cap R$ is a principal ideal $\idr{g}$ in $R$
generated by an irreducible polynomial $g$ in $R$.

\begin{proof}[Proof of Prop.\ \ref{prop:contractions-are-principal}]
      First suppose that $f\in S$ is irreducible. As we just showed, there is an irreducible
      $g\in R$ such that $\ids{f}\cap R=\idr{g}$. Furthermore, if $n\ge1$ then $\ids{f^n}$ is
      $\ids{f}$-primary, and so $\ids{f^n}\cap R$ is $\idr{g}$-primary, hence equals
      $\idr{g^k}$ for some $k\ge1$.

      The result is obvious if $f=0$, so suppose that $0\ne f\in S$, and let $f=f_1^{n_1}\cdots
      f_r^{n_r}$ be its factorization in $S$ into powers of distinct irreducibles $f_j$. Then
      there are irreducible polynomials $g_j\in R$ and $k_j\ge1$ such that $\ids{f_j^{n_j}}\cap
      R=\idr{g_j^{k_j}}$. Hence
      \begin{align*}\label{eqn:contraction}
        \ids{f}\cap R &= \ids{f_1^{n_1}\cdots f_r^{n_r}}\cap R = \bigl(\ids{f_1^{n_1}}\cap
                        \cdots \ids{f_r^{n_r}} \bigr)\cap R  \\
        &= \bigl( \ids{f_1^{n_1}} \cap R \bigr)\cap \dots \cap\bigl(\ids{f_r^{n^{}_r}}\cap R\bigr) 
          \\
        &= \idr{g_1^{k_1}}\cap\dots\cap\idr{g_r^{k_r}} = \idr{\operatorname{LCM}( g_1^{k_1},\dots,
         g_r^{k_r})}, 
      \end{align*}
      proving that $\ids{f}\cap R$ is principal.
\end{proof}

\begin{remarks}\label{rem:prim}
   (1) It is possible for distinct principal prime ideals in $S$ to contract to the same prime
   ideal in $R$. As a simple example, let $d=1$, $N=2$, $f_1(x)=x^2-x-1$, and
   $f_2(x)=x^2+x-1$. Then each is irreducible in $S$, but both $\ids{f_1}$ and $\ids{f_2}$
   contract in $R=\ZZ[2\ZZ]$ to $\idr{x^4-3x^2+1}$, where $x^4-3x^2+1$ is irreducible in $\ZZ[2\ZZ]$ (but of
   course not in $\ZZ[\ZZ]$). In the proof this is accounted for by using the least common
   multiple $\operatorname{LCM}$ in the last line of the displayed equation above.

   (2) A polynomial is \emph{primitive} if the greatest common divisor of its coefficients is
   1. If $0\ne f\in S$ is a nonconstant primitive polynomial with factorization $f=f_1^{n_1}\cdots f_r^{n_r}$ into
   powers of distinct irreducible polynomials, then by Gauss's Lemma each $f_j$ is primitive as
   well. Furthermore, $\ids{f_j}\cap R = \idr{g_j}$, where each $g_j$ is nonconstant and
   primitive. It then follows from the proof that $\ids{f}\cap R$ is generated by a primitive
   element of $R$.

   (3) There is a completely different proof of Prop.\ \ref{prop:contractions-are-principal}
   using entropy that is valid for all polynomials in $S$ except for those of a very special
   and easily determined form. Recall that the entropy of $\af$ is the logarithmic Mahler
   measure $\mahler(f)$ defined in 
   \eqref{eqn:mahler-measure}. A \emph{generalized cyclotomic polynomial} in $S$ is one of the
   form $\bx^{\bn}c(\bx^{\bk})$, where $c$ is a cyclotomic polynomial in one variable and
   $\bk\ne\boldsymbol{0}$. Smyth \cite{SmythKronecker} proved that $\mahler(f)=0$ if and only
   if $f$ is, up to sign, a product generalized cyclotomic polynomials. Assume that $f\in S$ is
   not such a polynomial, so that the entropy of $\af$ is strictly positive. A simple argument
   using cosets of $\nzd$ shows that $r_N(\xf)$ also has positive entropy. Now
   $r_N(\xf)=X_{\fa_N}$ by Lemma \ref{lem:contracted-ideal}, where $\fa_N=\ids{f}\cap R$. But
   an ideal $\fa$ in $R$ for which the shift action of $\nzd$ on $X_{\fa}$ has positive entropy must be
   principal \cite{LSW}*{Thm.\ 4.2}.
\end{remarks}

\section{Absolutely irreducible factorizations and Gauss's Lemma}\label{sec:factorizations}

Suppose that $f\in\zzd$ is nonconstant and irreducible. Its factorization into absolutely irreducible polynomials in an extension field of $\QQ$ will play a decisive role. A generalization of Gauss's Lemma to number fields enables us to deal with the algebraic properties of the coefficients of the factors.

Two polynomials in $\czd$ are \emph{distinct} if one is not a nonzero scalar multiple of the other. An element $\phi\in\czd$ is \emph{adjusted} if $\bzero$ is an extreme point of its Newton polytope $\SN_\phi$, and is \emph{monic} if it is both adjusted and $\widehat{\phi}(\bzero)=1$.

A polynomial in $\czd$ is \emph{absolutely irreducible} if it is irreducible in the unique factorization domain $\czd$. Hence every non-unit $f\in\czd$ has some factorization $f=\phi_1\cdots\phi_r$ into absolutely irreducible factors $\phi_j$. The method of Galois descent ~\cite{Conrad} shows that, after multiplying the factors by suitable constants, there is a finite normal extension $\KK$ of $\QQ$ such that each $\phi_j\in\kzd$, and also that the coefficients of the $\phi_j$ generate $\KK$, so that $\KK$ is the splitting field of $f$.
Furthermore an elementary argument shows that if $f$ is adjusted, then we can multiply the $\phi_j$ by units in $\kzd$ so that each $\phi_j$ is monic, $\SN_{\phi_j}\subset\nf$, and $f=\fh(\bzero)\phi_1\cdots\phi_r$.

\begin{remarks}
	(1) When $d=1$ this factorization is into the linear factors guaranteed by the fundamental theorem of algebra.
	
	(2) A simple sufficient condition for $\phi$ to be absolutely irreducible is that $\SN_\phi$ is not the nontrivial Minkowski sum of two integer polytopes (see \cite{Gao} for applications of this idea).
	
	(3) There are reasonably good factoring algorithms which, on input $f$, produce a monic irreducible polynomial in $\ZZ[x]$ with root $\theta$ and an absolutely irreducible $\phi\in \QQ(\theta)[\zd]$ such that $f =\sigma_1(\phi)\sigma_2(\phi)\cdots\sigma_r(\phi)$, where the $\sigma_j$ are all the distinct field embeddings of $\QQ(\theta)$ into $\CC$ (see \cite{Duval} for an overview of these methods). 
\end{remarks}

The following shows that, unlike factoring, divisibility is not affected when passing to an extension field.

\begin{lemma}\label{lem:factoring}
   Suppose that $\LL$ is an extension of the field $\KK$ and that $f,g\in\kzd$. Then $f$ divides $g$ in $\kzd$ if and only if $f$ divides $g$ in $\LL[\zd]$.
\end{lemma}

\begin{proof}
   For the nontrivial direction, suppose there is an $h\in\LL[\zd]$ such that $fh=g$. Equating coefficients of like monomials gives a system of $\KK$-linear equations in the coefficients of $h$. Since this system has a solution over $\LL$, Gaussian elimination shows that this (unique) solution is actually over $\KK$, and so $h\in\kzd$.
\end{proof}

\begin{proposition}\label{prop:absfactors}
   Let $f\in\zzd$ be nonconstant, adjusted, and irreducible in $\zzd$. Then there is a finite normal extension field $\KK$ of $\QQ$ and monic absolutely irreducible polynomials $\phi_1,\dots,\phi_r\in\kzd$ such that $f = \fh(\bzero)\phi_1\cdots\phi_r$ and $\SN_{\phi_f}\subset\nf$ for $1\le j\le r$. Furthermore, the Galois group $\Gal(\KK:\QQ)$ acts transitively on the set of factors $\phi_j$, and these factors are pairwise distinct.
\end{proposition}

\begin{proof}
   Our earlier discussion shows there is a factorization $f = \fh(\bzero)\phi_1\cdots\phi_r$ over the splitting field $\KK$ of $f$, where each $\phi_j$ is monic and $\SN_{\phi_j}\subset \nf$ for $1\le j\le r$. Suppose that $\sigma\in \Gal(\KK:\QQ)$. Since $\sigma(f)=f$, it follows that $\sigma$ must permute the absolutely irreducible factors up to multiplication by units. But if $\sigma(\phi_j) = c \bx^\bn \phi_k$, then $\bn=\bzero$ since the factors are adjusted and $c=1$ since they are monic. Hence $\sigma$ permutes the factors themselves. If there were a proper subset of factors that is invariant under $\Gal(\KK:\QQ)$, then their product $\psi$ would be in $\QQ[\zd]$ since its coefficients are invariant under $\Gal(\KK:\QQ)$. But then $\psi$ would be a proper divisor of $f$ in $\QQ[\zd]$ by Lemma \ref{lem:factoring}, contradicting irreducibility of $f$ by Gauss's Lemma. A similar argument shows that each factor appears with multiplicity one.
\end{proof}

We now give a brief sketch of the extension of Gauss's Lemma to number fields and the consequences we use. Let $\KK$ be a finite extension of $\QQ$, and $\OK$ be the ring of algebraic integers in $\KK$. A \emph{fractional ideal} $\fa$ in $\KK$ is a nonzero $\OK$-submodule such that there is an integer $b$ for which $b\fa\subset\OK$. Fractional ideals can be added and multiplied, with $\OK$ being the multiplicative identity. A fractional ideal contained in $\OK$ is an ideal in the usual ring-theoretic sense. The pivotal result is that the set of fractional ideals form a group, the set of principal fractional ideals (those of the form $\OK \beta$ for some $\beta\in\KK$) form a subgroup, and the quotient of these groups is a finite abelian group called the class group which measures how far $\OK$ is from being a principal ideal domain.

Let $\phi\in\kzd$. Define the \emph{content} $\ck(\phi)$ to be the factional ideal in $\KK$
generated by the coefficients of $\phi$. Say that $\phi$ is \emph{primitive} if
$\ck(\phi)=\OK$. It is easy to check that although content depends on the ambient field $\KK$,
primitivity does not: if $\phi\in\kzd$ and $\phi\in\lzd$, then $\ck(\phi)=\OK$ if and only if
$\cl(\phi)=\OL$. For a proof of the following version of Gauss's Lemma, see \cite{Cassels}*{Chap.\ 6, Lemma 2.1},  or the more leisurely account in \cite{MagidinMcKinnon}*{Thm.\ 8.2}.

\begin{theorem}[Gauss's Lemma for number fields]\label{thm:gauss-lemma}
   Let $\KK$ be a number field and $\phi,\psi\in\kzd$. Then
   $\ck(\phi\psi)=\ck(\phi)\ck(\psi)$. In particular, if $\phi,\psi\in\okzd$ then $\phi\psi$
   is primitive if and only if both $\phi$ and $\psi$ are primitive. If $\phi,\psi\in\okzd$ are
   primitive, and if $\phi=\beta\psi$ for some $\beta\in\KK$, then $\beta$ is a unit in $\OK$. 
\end{theorem}

\begin{remark}\label{rem:fn-is-primitive}
   Suppose that  $f\in\zzd$ is primitive and that $N\ge1$. Let $\zeta_N=e^{2\pi i/N}$, which is a unit in $\QQ(\zn)$. Hence each rotate $\bo\cdot f$, where $\bo\in\Omnd$, is primitive in $\QQ(\zeta_N)[\zd]$. The preceding theorem then shows that the product $\fn$ of these rotates is also primitive in $\QQ(\zn)[\zd]$, and hence in $\zzd$ (since primitivity is independent of ambient field), a fact we already observed in Remark \ref{rem:prim}(2).
\end{remark}

\section{Decimated polynomials and decimated actions}\label{sec:decimated-actions}

Let $f\in\zzd$ be irreducible. Here we explain the relationship between the $N$th decimation $\fn$ of $f$ and the generator $g_N$ of the contracted ideal $\id{f}\cap\znzd$ that defines the $N$th decimation $r_N(X_f)$ of $(X_f,\af)$.  Roughly speaking, $g_N$ is a constant times the product of all distinct rotates by elements of $\Omnd$ of the absolutely irreducible factors $\phi_j$ of $f$ as described in Proposition \ref{prop:absfactors}. Each rotate appears with the same multiplicity $e_N$ that can be computed from the $\phi_j$. Thus $\fln=c\, g_N^{e_N}$, and an application of Gauss's Lemma shows that we may take $c=1$. Furthermore, there is an integer $Q(f)$,  that can also be computed from the $\phi_j$, such that $\fn=g_N^{}$ for all $N$ relatively prime to $Q(f)$. Examples will illustrate the two sources of the multiplicity $e_N$.

In what follows we let $\zn=e^{2\pi i/N}$, which is a generator of $\Omn$.

\begin{lemma}\label{lem:integral}
   If $f\in\zzd$ then $\fn\in\znzd$.
\end{lemma}

\begin{proof}
   Since $\fn=\prod_{\bo\in\Omnd} \bo\cdot f$, it follows that $\fn=\bo\cdot\fn$ for every $\bo\in\Omnd$. Suppose that $\fnh(\bk)\neq 0$. Then since
   \begin{displaymath}
      \fnh(\bk)=(\bo\cdot\fn)\widehat{\ }\,(\bk)=\bo^{\bk}\,\fnh(\bk),
   \end{displaymath}
   we see that $\bo^{\bk}=1$ for every $\bo\in\Omnd$, and hence $\bk\in N\ZZ^d$. Thus $\fn\in\QQ(\zn)[N\ZZ^d]$.
   
   The Galois group $G := \Gal\bigl(\QQ(\zn): \QQ\bigr)$ acts on $\Omnd$ coordinate-wise. If $\sigma\in G$, then $\sigma(\bo\cdot f)=\sigma(\bo)\cdot f$ since $f$ has integer coefficients. Thus $\sigma$ permutes the rotates of $f$, and so $\sigma(\fn)=\fn$ for every $\sigma\in G$. It follows that the coefficients of $\fn$ are both rational and algebraic integers, and so $\fn\in\znzd$.
\end{proof}

\begin{lemma}\label{lem:divides}
   Let $f\in\zzd$ and $g_N$ be a generator of the contracted ideal $\id{f}\cap\znzd$. Then $g_N$ divides $\fn$ in $\znzd$.
\end{lemma}

\begin{proof}
   Since $f$ is one of the factors in forming $\fn$, it follows that $f$ divides $\fn$ in $\QQ(\zn)[\zd]$. Hence $f$ divides $\fn$ in $\QQ[\zd]$ by Lemma \ref{lem:factoring}. The coefficients of $\fn/f$ are both rational and algebraic integers, and so $\fn/f \in\zzd$. Hence $\fn\in\id{f}\cap\znzd$, and it is thus divisible by the generator $g_N$.
\end{proof}

\begin{remark}\label{rem:pos}
	Since the generator of a principal ideal is unique only up to units, it will be convenient to have a convention to pick a generator. In what follows we will assume that $f$ is adjusted and that $\fh(\bzero)>0$. Then clearly $\fn$ has the same properties. By the previous lemma, we can also assume that $g_N$ is adjusted, that $\SN_{g_N}\subset\SN_{\fn}$, and that $\gnh(\bzero)>0$.
\end{remark}

Before continuing, we remark that if $f$ is a constant integer $n$, then $\fn=n^{N^d}$ while $g_N=n$. Let us call a polynomial $f\in\zzd$ \emph{nonconstant} if $|\supp f|>1$, and it is these we now turn to.

Let $f\in\zzd$ be adjusted. Define its \emph{support group} $\Gf$ to be the subgroup of $\zd$ generated by $\supp f$. It is easy to check that the support group is independent of which extreme point of $\nf$ is used to adjust $f$. We say that $f$ is \emph{full} if $\Gf=\zd$.

The following shows that in some cases, including $f(x,y)=1+x+y$ from Example ~\ref{exam:ledrappier}, $\fn=g_N^{}$ for all $N\ge1$.

\begin{proposition}\label{prop:absfull}
   Let $f\in\zzd$ be adjusted, irreducible, and full. Further assume that $f$ is absolutely irreducible in $\czd$. Then $\fn =g_N^{}$ for every $N\ge1$.
\end{proposition}

\begin{proof}
   Since the map $f\mapsto\bo\cdot f$ is a ring isomorphism of $\czd$, each $\bo\cdot f$ is absolutely irreducible. Suppose that $\bo\cdot f=\bo'\cdot f$. Since $\bo^{\bk}\fh(\bk)=(\bo')^{\bk}\fh(\bk)$, it follow that $\bo^{\bk}=(\bo')^{\bk}$ for all $\bk\in\supp f$, hence for all $\bk\in\Gf=\zd$, and so $\bo=\bo'$. Thus the rotates of $\bo\cdot f$ for $\bo\in\Omnd$ are pairwise distinct absolutely irreducible polynomials in $\czd$ whose product is $\fn$.
   
   By Lemma \ref{lem:divides}, $\gn$ divides $\fn$ in $\czd$. Hence some rotate $\bo\cdot f$ divides $\gn$. Since $\gn\in\znzd$, it is invariant under all rotations in $\Omnd$. Hence $\gn$ is divisible by all rotates $\bo\cdot f$, and so $\gn$ and $\fn$ have the same absolute factorizations in $\czd$, and hence $\fn=c\,\gn$ for some constant $c\in\CC$. Recalling our conventions in Remark \ref{rem:pos}, comparing constant terms shows that $c= \fh(\bzero)^{N^d}/\,\gnh(\bzero)\in\QQ$. But $\fn$ and $\gn$ are both primitive in $\znzd$, and so $c=\pm1$, and our convention on positivity of constant terms then gives $c=1$.
\end{proof}

The following example shows that when the polynomial is not full there can be multiplicity $\en>1$.

\begin{example}\label{exam:1+x+y2}
   Let $d=2$ and $f(x,y)=1+x+y^2$. Since $\nf$ is not a nontrivial Minkowski sum of integer polytopes, we see that $f$ is absolutely irreducible. Suppose that $N$ is odd. Since $-1\notin\Omn$, all rotates $\bo\cdot f$ for $\bo\in\Omn^2$ are distinct, and the same arguments as in the previous proposition show that $\fn=\gn^{}$.
   
   However, if $N$ is even, then $-1\in\Omn$ and the rotate of $f$ by $(\omega_1,\omega_2)$ equals that by $(\omega_1,-\omega_2)$. As we will see in Proposition \ref{prop:support}, the product of the distinct rotates of $f$ equals $\gn$, and so $\fn=\gn^2$ when $N$ is even. 
\end{example}

Next we characterize when rotates can coincide.

\begin{lemma}\label{lem:stabilizer}
   Let $\phi\in\czd$ be adjusted, and $\Gp$ be its support group. Then the dual of the stabilizer group $\sn(\phi) := \{\bo\in\Omnd:\bo\cdot\phi=\phi\}$ is $\zd/(\Gp+\nzd)$. Two rotates of $\phi$ differ by a multiplicative unit in $\czd$ if and only if they are equal. If $\Gp$ has finite index $K$ in $\zd$, then $\sn(\phi)$ is trivial for every $N$ relatively prime to ~$K$.
\end{lemma}

\begin{proof}
   Suppose that $\bo\in\sn(\phi)$. Since $\ph(\bk)=(\bo\cdot\phi)\widehat{\ }\,(\bk)=\bo^\bk\ph(\bk)$, it follows that $\bo^\bk=1$ for every $\bk\in\supp \phi$. Hence $\bo$ annihilates $\Gp$ as well as $\nzd$, thus their sum. Conversely, every $\bo$ annihilating $\Gp+\nzd$ must be in $\sn(\phi)$. Hence the annihilator of $\sn(\phi)$ equals $\Gp+\nzd$, and so its dual group is $\zd/(\Gp+\nzd)$.
   
   The multiplicative units in $\czd$ have the form $c\,\bx^{\bn}$ for some $c\in\CC$, so the second statement is obvious since $\phi$ is adjusted.
   
   Suppose that $\Gp$ has finite index $K$ in $\zd$. If $N$ is relatively prime to $K$, then multiplication by $N$ on $\zd/\Gp$ is injective, hence surjective. Thus modulo $\Gp$ every element in $\zd$ is a multiple of $N$, and hence $\Gp+\nzd=\zd$.
\end{proof}

\begin{proposition}\label{prop:support}
   Let $f\in\zzd$ be adjusted and irreducible, and further assume that $f$ is absolutely irreducible in $\czd$. Then $\fln=\gn^{\en}$, where $\en=|\sn(f)|=|\zd/(\Gf+\nzd)|$.
\end{proposition}

\begin{proof}
   Recall our conventions in Remark \ref{rem:pos}. Since $\gn$ divides $\fn$, it must be divisible by at least one (absolutely irreducible) rotate of $f$. Invariance of $\gn$ by every rotate in $\Omnd$ shows that $\gn$ is therefore divisible by the product $h$ of all the distinct rotates of $f$. The arguments in Lemmas \ref{lem:integral} and \ref{lem:divides} apply to show that $h\in\id{f}\cap\znzd$. Thus $\gn$ divides $h$ in $\czd$ as well, and so $\gn=c\,h$ for some $c\in\CC$. Evaluating constant terms shows that $c\in\QQ$. Since $\gn$ is irreducible in $\znzd$, it is primitive. Each rotate of $f$ is primitive in $\QQ(\zn)[\zd]$, and so $h$ is primitive by Theorem \ref{thm:gauss-lemma}. Hence $c=\pm1$, and then $c=1$ follows from our sign conventions. By Lemma~ \ref{lem:stabilizer}, each rotate of $f$ is repeated exactly $\en$ times, and so $\fln=\gn^{\en}$.
\end{proof}

When $f$ is absolutely irreducible, the only source of multiplicity $\en>1$ is its support group. However, if $f$ has several absolutely irreducible factors, a new source of multiplicity can occur, namely that one factor could rotate to another factor. This possibility is illustrated in the following three examples.

\begin{example}\label{exam:1-2x2}
   Let $d=1$ and 
   \begin{displaymath}
      f(x)=1-2x^2=(1+\sqrt{2}x)(1-\sqrt{2}x)=\phi_1(x)\phi_2(x).
   \end{displaymath}
   Let $\sigma\in\Gal(\QQ(\sqrt{2}):\QQ)$ be given by $\sigma(\sqrt{2})=-\sqrt{2}$. Then $\sigma(\phi_1)=\phi_2=(-1)\cdot\phi_1$. Now $\fn$ is the product of $\zn^j\cdot\phi_k$ for $0\le j<N$ and $k=1,2$. If $N$ is odd, then $-1\notin\Omn$ and so all $2N$ factors are distinct. Our earlier arguments then show that $\fn=\gn$. However, if $N$ is even, then $-1\in\Omn$, and the set of rotates of $\phi_1$ coincide with set of those of $\phi_2$, and so $\fn=\gn^2$ for even $N$. Here $f$ is an irreducible polynomial with a pair of roots whose ratio is a nontrivial root of unity.
\end{example}

The commingling of absolutely irreducible factors under rotations can happen in more subtle ways.

\begin{example}\label{exam:4th-deg}
   Let $d=1$ and $f(x)=1-2x+4x^2-3x^3+x^4$, which is full and irreducible in $\ZZ[\ZZ]$. Let $\lambda=(1+\sqrt{5})/2$, $\mu=(1-\sqrt{5})/2$, and $\zeta=\zeta_5$. The absolutely irreducible factorization of $f$ is
   \begin{displaymath}
         f(x)=(1-\zeta\lambda x)((1-\zeta^4 \lambda x)(1-\zeta^2\mu x)(1-\zeta^3 \mu x)
             = \phi_1(x)\phi_2(x)\phi_3(x)\phi_4(x).
   \end{displaymath}
   Note that $\zeta^3\cdot\phi_1=\phi_2$ and that $\zeta\cdot\phi_3=\phi_4$. If $N$ is relatively prime to 5, then $\zeta\notin\Omn$, and so all $4N$ rotates are distinct and $\fn=\gn$ as before. However, if $5\mid N$ then $\zeta\in\Omn$ and each rotate is repeated twice, and so $\fn=\gn^2$ in this case.
   
   What is driving this example is the inclusion $\QQ(\sqrt{5})\subset\QQ(\zeta)$, and so the Galois automorphism $\sqrt{5}\mapsto -\sqrt{5}$ of $\QQ(\sqrt{5})$ is the restriction of the automorphism $\zeta\mapsto\zeta^2$ of $\QQ(\zeta)$.
\end{example}

\begin{remark}\label{rem:degenerate}
   Irreducible polynomials in $\ZZ[x]$ having distinct roots whose ratio is a root of unity, such as those in the previous two examples, are called \emph{degenerate}. Such polynomials have an extensive literature (see for instance \cite{recurring}*{\S1.1.9}), and appear in the celebrated Skolem-Mahler-Lech Theorem that the set of indices at which a recurring sequence of integers vanishes is, modulo a finite set, the union of arithmetic progressions \cite{Berstel}.
   
   There is a simple way to detect whether $f(x)\in\ZZ[x]$ is degenerate. Introduce a new variable $t$, and compute the resultant $g(x)\in\ZZ[x]$ of the polynomials $f(tx)$ and $f(t)$ with respect to $t$, which can be done efficiently using rational arithmetic. The roots of $g(x)$ are the ratios of all pairs of roots of $f$. Thus $f(x)$ is degenerate if and only if $g(x)$ contains a nontrivial cyclotomic factor. Applying this to $f(x)$ from the previous example gives
   \begin{displaymath}
      g(x)=(x-1)^5(x^4-4x^3+6x^2+x+1)(x^4+x^3+6x^2-4x+1)(x^4+x^3+x^2+x+1).
   \end{displaymath}
   The last factor reveals that $f(x)$ has two roots whose ratio is a nontrivial 5th root of unity.
\end{remark}

\begin{example}\label{exam:twisted-ledrappier}
   Let $d=2$ and $f(x,y)=1-x-y-xy+x^2+y^2$, which is full and irreducible in $\ZZ[\ZZ^2]$. Let $\zeta=\zeta_3$. The absolutely irreducible factorization of $f$ is
   \begin{displaymath}
      f(x,y) = (1+\zeta x + \zeta^2 y)(1+\zeta^2 x + \zeta y)=\phi_1(x,y)\phi_2(x,y).
   \end{displaymath}
   Here $\phi_1$ is mapped to $\phi_2$ by the element $\sigma$ in $\Gal(\QQ(\zeta):\QQ)$ mapping $\zeta$ to $\zeta^2$, and also $\sigma(\phi_1)=\phi_2=(\zeta,\zeta^2)\cdot\phi_1$. By the now familiar arguments, if $N$ is relatively prime to 3 then $\zeta\notin\Omn$, and so all rotates are distinct and hence $\fn=\gn$. However, if 3 divides $N$, then distinct rotates are repeated twice, and so $\fn=\gn^2$. For instance
   \begin{displaymath}
      f_{\<3\>}= (1+3x^3+3y^3+3x^6-21x^3y^3+3y^3+x^9+3x^3y^6+3x^6y^3+y^9)^2 = g_3^2.
   \end{displaymath}
\end{example}

With these examples in mind, we come to the main result of this section.

\begin{theorem}\label{thm:contraction}
   Let $f\in\zzd$ be irreducible, which we may assume is adjusted with positive constant term. For every $N\ge1$ there is an irreducible $\gn\in\znzd$ and $\en\ge1$ such that
   \begin{displaymath}
      \idzd{f}\cap \znzd = \idnzd{\gn} \text{\quad and\quad} \fln=\gn^{\en}.
   \end{displaymath}
   The multiplicity $\en$ can be computed from the absolutely irreducible factorization of ~$f$ in $\czd$. If the support of $f$ generates a finite-index subgroup of $\zd$, then there is an integer $Q(f)$, which can also be computed from the absolutely irreducible factors of~ $f$, such that $\en=1$ for every $N$ that is relatively prime to $Q(f)$. Finally,
   \begin{displaymath}
      \idzd{f^k}\cap\znzd = \idnzd{\gn^k}
   \end{displaymath}
   for every $k\ge1$.
\end{theorem}

\begin{proof}
   Recall our conventions in Remark \ref{rem:pos}. Let $\KK$ be the splitting field of $f$, and $f=\fh(\bzero)\phir$ be the factorization of $f$ using monic absolutely irreducible $\phi_j\in\kzd$ from Proposition \ref{prop:absfactors}. Let $\Phi=\{\phi_1,\dots,\phi_r\}$. Since the $\phi_j$ are monic, $\Gal(\KK:\QQ)$ permutes the elements of $\Phi$, and this action is transitive by irreducibility of $f$.
   
   Now fix $N\ge1$. Then $\KK(\zn)$ is a normal extension of $\QQ$. Let $G=\Gal(\KK(\zn):\QQ)$. Consider the set $\Omnd\times \Phi$. The group $\Omnd$ acts on this set via $\bo'\cdot (\bo,\phi_j)=(\bo'\bo,\phi_j)$. The group $G$ also acts on this set via $\sigma\cdot(\bo,\phi_j)=(\sigma(\bo),\sigma(\phi_j))$. More precisely, $\sigma\in G$ acts of the first coordinate using its restriction to $\QQ(\zn)$ and on the second coordinate using its restriction to $\KK$. These actions combine to give an action of the semidirect product $G\ltimes \Omnd$ defined using the action of $G$ on $\Omnd$, so that $\sigma\bo=\sigma(\bo)\sigma$.
   
   Define an equivalence relation $\sim$ on $\omndp$ by $(\bo,\phi_j)\sim(\bo',\phi_k)$ if and only if $\bo\cdot \phi_j=\bo'\cdot\phi_k$. It is routine to verify that $\gond$ preserves equivalence classes. Since $\Gal(\KK:\QQ)$ acts transitively on $\Phi$, it follow that $\gond$ acts transitively on $\omndp$. Hence all equivalence classes have the same cardinality, say $\en\ge1$. Pick one representative $(\bo,\phi_j)$ from each equivalence class, and let $\gtn$ be the product of the corresponding polynomials $\bo\cdot\phi_j$.
   
   Observe that by its construction $\gtn$ is invariant under $\gond$. Invariance under $\Omnd$ implies that $\gtn\in\KK(\zn)[\nzd]$, and invariance under $G$ further implies that $\gtn\in\QQ[\nzd]$. Then transitivity of $\gond$ on $\omndp$ shows that $\gtn$ is irreducible in $\QQ[\nzd]$.
   
   We have that $\fn=\fh(\bzero)^{N^d}\,\gtn^{\,\en}$. Let $q$ be the least positive integer such that $q \gtn\in\znzd$, so that $\gn:=q\gtn$ is primitive. Then 
   \begin{displaymath}
      \fn=\bigl( \fh(\bzero)^{N^d} / q^{\en}\bigr) \gn^{\en}.
   \end{displaymath}
   But both $\fn$ and $\gn^{\en}$ are primitive with positive constant terms, and hence $\fn=\gn^{\en}$.
   
   We now turn to computing $\en$. Each of the absolutely irreducible factors $\phi_j$ has the same support since they are all Galois conjugates. Let $\Gp$ denote the common support group of each. By Lemma \ref{lem:stabilizer}, each contributes multiplicity $|\zd/(\Gp+\nzd)|$. Further multiplicity arises if one factor can be rotated by an element of $\Omnd$ to another. This property divides $\Phi$ into equivalence classes, with all classes having the same cardinality $s$. It then follows that $\en= |\zd/(\Gp+\nzd)|s$.
   
   Next, we determine sufficient conditions on $N$ so that $\en=1$. Assume that $\Gf$ has finite index in $\zd$. Clearly $\Gf\subset\Gp$, and so $\Gp$ also has finite index. By Lemma ~\ref{lem:stabilizer}, if $N$ is relatively prime to the index $[\zd:\Gp]$ of $\Gp$, then $|\zd/(\Gp+\nzd)|=1$.
   
   To analyze when one $\phi_j$ can rotate to another, we need to consider the group $\Omega_\KK$ of roots of unity in the splitting field $\KK$ of $f$. This is a finite cyclic group, and so equals $\Omega_n$ for some $n\ge1$. Now $[\QQ(\zeta_n):\QQ]=\varphi(n)$, where $\varphi$ denotes the Euler function. Since $\QQ(\zeta_n)\subset \KK$, it follows that  $\varphi(n)\le [\KK:\QQ]$. A simple argument shows that $\varphi(n)\ge\sqrt{n}/2$ for all $n\ge 1$, and so $n\le 4[\KK:\QQ]^2$. Hence if $N$ is relatively prime to $(4[\KK:\QQ]^2)!$, then $\Omn\cap\Omega_\KK=\{1\}$. For such an $N$ suppose that $\bo\cdot\phi_i=\phi_j$ for some $\bo\in\Omnd$. For each $\bk\in\supp\phi_i=\supp\phi_j$ we have that $\bo^{\bk}\widehat{\phi}_i(\bk)=\widehat{\phi}_j(\bk)$, and so
   \begin{displaymath}
       \bo^{\bk}= \widehat{\phi}_j(\bk)/\widehat{\phi}_i(\bk)\in\Omn\cap\Omega_\KK = \{1\}.
   \end{displaymath}
   But this implies that $\phi_i=\phi_j$.
   
   Putting these together, we let $Q(f)=[\zd:\Gp](4[\KK:\QQ]^2)!$, and conclude that if $N$ is relatively prime to $Q(f)$ then $\en=1$. 
\end{proof}

\section{Remarks and Questions}\label{sec:problems} 

Here we make some further remarks and ask several questions related to decimations.

\subsection{More general lattices} \label{subsec:general}
Let us call a finite-index subgroup of $\zd$ a \emph{lattice}. We have used the sequence $\{\nzd\}$ of lattices to define decimation, but these definitions easily extend to all lattices. Let $\Lam\in\zd$ be a lattice, and let $\Omega_{\Lam}$ denote the dual group of $\zd/\Lam$, which has cardinality $[\zd\colon\Lam]$, the index of $\Lam$ in $\zd$. Define $f_{\langle\Lam\rangle}=\prod_{\bo\in\Omega_{\Lam}}\bo\cdot f$, and
\begin{equation}\label{eqn:partialdecimation}
   \ssL_{\Lam}f = \ssE_{[\zd\colon\Lam]}\Bigl( \frac{1}{[\zd\colon\Lam]}\log|\fh_{\langle\Lam\rangle}|\Bigr).
\end{equation}
For a sequence $\{\Lam_N\}$ of lattices, let us say $\Lam_N\to\infty$ if for every $r>0$ we have that $\{\bn\in\Lam_N:\|\bn\|<r\}=\{\bzero\}$ for all large enough $N$.

If we replace the ``square'' lattices $\nzd$ with ``rectangular'' lattices of the form
\begin{displaymath}
   \Lam_N = a_{N}^{(1)}\ZZ \oplus \dots \oplus a_{N}^{(d)}\ZZ,
\end{displaymath}
then a straightforward modification of our proof shows that $\ch(\ssL_{\Lam_N} f)\to\df$ uniformly on $\nf$ provided that $a_N^{(k)}\to\infty$ as $N\to\infty$ for each $1\le k \le d$.

However, the analogous question for general lattices is much more subtle, since the rescaling argument that is basic to our proof has no obvious extension. Nevertheless, in recent unpublished work Hanfeng Li has been able to establish the uniform convergence of $\ch(\ssL_{\Lam_N} f)$ to $\df$ for every sequence $\{\Lam_N\}$ of lattices with $\Lam_N\to\infty$.

We can also investigate decimations by lattices with different limiting behavior. Let $\SC(\zd)$ denote the set of subgroups of $\zd$. We can give a topology to $\SC(\zd)$ by declaring two subgroups to be close if they agree on a large ball around $\bzero$. For example, in this topology $\Lam_N\to \{\bzero\}$ means that $\Lam_N\to\infty$ as above.  This is a special case of the Chabauty topology on the set $\SC(G)$ of closed subgroups of a locally compact group $G$. This topology is named after Claude Chabauty, who in 1950 introduced it \cite{Chab} to generalize Mahler's compactness criterion \cite{MahlerCompactness} for lattices in $\RR^d$ to lattices in locally compact groups. The Chabauty space $\SC(G)$ has been investigated by many authors, for instance by Cornulier \cite{Corn} when $G$ is abelian. Even for familiar groups their Chabauty space can be intricate to analyze. For example, Hubbard and Pourezza \cite{HP} used a tricky argument to prove that $\SC(\RR^2)$ is homeomorphic to the four-dimensional sphere.

\begin{question}\label{ques:general-lattices}
   Let $\Delta$ be a subgroup of $\zd$ of infinite index. Suppose that $\{\Lam_N\}$ is a sequence of lattices such that $\Lam_N\to\Delta$ as $N\to\infty$. Do the functions $\ch(\ssL_{\Lam_N} f)$ always converge on $\nf$, and if so what is the limit function in terms of $f$ and $\Delta$?
\end{question}

\subsection{Exponential size of decimation coefficients}

In Example \ref{exam:transcendental} we saw that if $f\in\CC[\ZZ]$ is allowed to have complex coefficients, then some of the coefficients of $\fn$ may have exponential size drastically different from that predicted by $\df$. However, if $f\in\ZZ[\ZZ]$ is restricted to have integer coefficients, then this behavior cannot happen, as indicated by Example \ref{exam:quasihyperbolic}. More precisely, using the diophantine results of Gelfond mentioned there, one can show that if $f\in\ZZ[\ZZ]$ has $\supp f = \{0,1,\ldots,r\}$ and $\varepsilon>0$, then for all sufficiently large $N$ we have that $|\fnh(kN)|$ is between $e^{N( \df(k)\pm\varepsilon)}$ for each $0\le k \le r$ for which $\fnh(kN)\ne 0$.

This raises the intriguing question of whether this extends to $f\in\zzd$ for \mbox{$d\ge 2$}, i.e., do all nonzero coefficients of $\fn$ have the approximate exponential size predicted by $\df$. The following gives a precise quantitative formulation.

\begin{question}\label{que:coefficcients}
   Let $f\in\zzd$. Fix $\br_0\in\nf$, and let $\varepsilon>0$. Are there $\delta>0$ and $N_0\ge1$ such that if $N\ge N_0$ and $\br\in N^{-d}\zd\cap\nf$ with $\|\br-\br_0\|<\delta$, and if $\ssL_N f(\br)\neq -\infty$, then $|\ssL_N f(\br)-\df(\br)|<\varepsilon$? Can $\delta$ and $N_0$ be chosen uniformly for $\br_0\in\nf$?
\end{question}

Some evidence for a positive answer comes from polynomials in two variables related to dimer models, as discussed in Remark \ref{rem:1+x+y}. Using the additional machinery afforded by the physical interpretation of the related partition function and the resulting subadditivity, the exponential size of the coefficients can be shown to obey the estimates in the question. In particular, this applies to $f(x,y)=1+x+y$, although we do not know of any direct argument for this.

\subsection{Continuity of $\exp[\df]$ in the coefficients of $f$}

Start by fixing a cube $B_n=\{-n,\ldots,n\}^d\subset\zd$. We can identify a polynomial $f\in\czd$ whose support is in $B_n$ with its coefficient function $\fh\in\CC^{B_n}$. Boyd \cite{BoydUniform} showed that the function $\CC^{B_n}\to[0,\infty)$ given by $\fh\mapsto \Mahler(f)=\exp[\mahler(f)]$ is continuous in the coefficients of $f$.

Recalling that $\mahler(f)$ is the maximum value of $\df$, this suggests looking at $\exp[\df]$, which is a nonnegative upper semicontinuous function on $B_n$ (the discontinuities occur at the boundary of $\nf\subset B_n$). A function $\phi\colon B_n\to\RR$ is upper semicontinuous if and only if its subgraph $\{(\bu,t)\in B_n\times\RR:t\le \phi(\bu)\}$ is closed in $B_n\times\RR$. The space $\USC(B_n)$ of all upper semicontinuous functions on $B_n$ carries a natural topology by declaring two elements to be close if their subgraphs are close in the Hausdorff metric on closed subsets of $B_n\times\RR$ (see \cite{Beer} for details). 

\begin{question}
   Is the map $\fh\to\exp[\df]$ from $\CC^{B_n}$ to $\USC(B_n)$ continuous?
\end{question}

\subsection{Nonprincipal actions}

Decimation makes sense for every algebraic $\zd$-action (indeed for every algebraic action of a countable residually finite group). Suppose that $\mathfrak{a}$ is an ideal in $\zzd$, and let $X_{\mathfrak{a}}$ be the dual group of $\zzd/\mathfrak{a}$ as described in ~\S\ref{sec:algebraic-decimations} with its associated algebraic $\zd$-action $\alpha_{\fa}$. The commutative algebra there shows that the $N$th decimation $r_N(X_{\mathfrak{a}})$ is defined by the contracted ideal $\mathfrak{a}\cap\znzd$. However, there is no obvious replacement for $g_N$ to measure growth when $\mathfrak{a}$ is not principal,

\begin{question}
   If $\mathfrak{a}$ is a nonprincipal ideal in $\zzd$, are there objects related to the contractions $\mathfrak{a}\cap\znzd$ which can be normalized to converge to a limiting object?
\end{question}

If $\mathfrak{a}$ is not principal, then the $\zd$-shift action $\alpha_{\fa}$ on $X_{\mathfrak{a}}$ has zero entropy. However, by restricting the shift to iterates close to lower dimensional subspaces of $\RR^d$ the action can have positive entropy \cite{BoyleLind}*{\S6}. This suggests that Question \ref{ques:general-lattices} may be relevant here.

Examining concrete examples may shed some light on this question. These include the case of commuting toral automorphisms (see \cite{KKS}*{\S6} for many such examples), the $\ZZ^2$-action defined by multiplication by 2 and by 3 on $\TT$ (corresponding to $\fa=\<x-2,y-3\>$), and the so-called space helmet example \cite{EinsiedlerLind}*{Example~5.8} (corresponding to $\fa=\<1+x+y,z-2\>$).

An important example of a different character is due to Ledrappier \cite{Ledrappier}, which corresponds to the nonprincipal ideal $\<1+x+y,2\>\subset\ZZ[\ZZ^2]$. This example has zero entropy as a $\ZZ^2$-action, but strictly positive entropy along every 1-dimensional subspace of $\RR^2$ (see \cite{BoyleLind}*{Example 6.4} for the explicit description). Another curious feature of this example is decimation self-similarity. Because $(1+x+y)^{2^n}=1+x^{2^n}+y^{2^n}$ when taken mod ~2, the $2^n$th decimation of the example, when rescaled by ~$2^n$, is just the original action.

\subsection{Computing entropy using decimations}

In his elegant proof that Mahler measure is continuous in the coefficients of polynomials $f\in\czd$, Boyd \cite{BoydUniform} side-stepped delicate issues about logarithmic singularities of $\log |f|$ by instead computing decimations of $f$ along powers of 2 using the Graeffe root-squaring algorithm. This enabled him to compute $\Mahler(f)$ to any prescribed accuracy using only a finite number of arithmetic operations on the coefficients of ~$f$.

This suggests a new approach to computing the entropy of algebraic actions using decimations and contracted ideals.

To describe this approach, let $\fa$ be an ideal in $\zzd$, and $\alpha_{\fa}$ be its associated algebraic $\zd$-action. Define the \emph{length} of $g\in \zzd$ to be $L(g):=\sum_{\bk\in\zd} |\gh(\bk)|$. For a subset $\fb$ of $\zzd$ put $L(\fb):=\min\{L(g)\colon 0\ne g \in\fb\}$. By convention we define $L(\{0\}):=\infty$. We define the \emph{asymptotic length} $\lam(\fa)$ of the ideal $\fa$ by
\begin{displaymath}
   \lam(\fa) :=\limsup_{N\to\infty} \frac{1}{N^d} \log L(\fa \cap \znzd).
\end{displaymath}

The following shows that for some principal ideals the asymptotic length equals the entropy of the associated action.

\begin{proposition}\label{prop:asymptotic-length}
   Suppose that $f\in\zzd$ is absolutely irreducible, adjusted, and has support that generates $\zd$. Then $\lam\bigl(\<f\>\bigr) =\h(\af)$.
\end{proposition}

\begin{proof}
   First note that by Theorem \ref{thm:contraction}, $\<f\>\cap\znzd=\<\fn\>_{\znzd}$. Choose $B_0$ so that $\supp f\subset \{-B_0,\dots,B_0\}^d$, and put $B=2B_0+1$. By \eqref{eqn:Mahler-bound} applied to the case $\bu=\bzero$, we see that
   \begin{displaymath}
      \Ht(\fn) \le 2^{dN^{d-1}B} \Mahler(f)^{N^d},
   \end{displaymath}
   and hence
   \begin{displaymath}
     L(\fn)\le (N^d B)^d 2^{dN^{d-1}B} \Mahler(f)^{N^d}.
   \end{displaymath}
   Thus
   \begin{displaymath}
      \frac{1}{N^d} \log L(\<f\>\cap  \znzd) \le \frac{1}{N^d} \log L(\fn) 
      \le \frac{d\log(N^d B)+d N^{d-1} B \log 2}{N^d} +\mahler(f),
   \end{displaymath}
   and so letting $N\to\infty$ we conclude that $\lam(\< f\>)\le \mahler(f)=\h(\af)$.
   
   To prove the reverse inequality, let $g\fn$ be an arbitrary nonzero element in $\<f\>\cap\znzd$. Then $\Mahler(g)\ge1$. Furthermore, $\Mahler(g\fn)\le L(g\fn)$ by \eqref{eqn:length-inequality}. Hence
   \begin{displaymath}
      \Mahler(f)^{N^d}=\Mahler(\fn)\le \Mahler(g)\Mahler(\fn) = \Mahler(g\fn)\le L(g\fn).
   \end{displaymath}
   Thus
   \begin{displaymath}
      \h(\af)=\mahler(f) = \frac{1}{N^d} \log \Mahler(\fn)\le \log L(g\fn),
   \end{displaymath}
   and so 
   \begin{displaymath}
      \lam(\<f\>)\ge\liminf _{N\to\infty} \frac{1}{N^d} \log L(\< f\> \cap \znzd) \ge \h(\af).
      \qedhere
   \end{displaymath}
\end{proof}

Under the assumptions on $f$ in the above result there is actually convergence to ~$\h(\af)$. However, convergence can fail if the exponents $e_N$ in Theorem \ref{thm:contraction} are at least 2 infinitely often. To give a simple example, let $d=1$ and $f(x)=x^2-2$. If $N$ is odd then $g_N(x)= x^{2N}-2^N$, while if $N$ is even then $g_N(x)= x^{N}- 2^{N/2}$. Thus $(1/N)\log L(\<f\>\cap\znzd)$ converges to $\frac{1}{2}\log 2$ along even $N$ and to $\log 2$ along odd ~$N$.

More seriously, Proposition \ref{prop:asymptotic-length} can fail when the support of $f$ does not generate a finite-index subgroup of $\zd$. Take for instance $f=2$, so that $\lam(\<f\>)=0$ but $\h(\af)=\log 2$.

\begin{question}\label{ques:entropy}
	What are necessary and sufficient conditions on $f\in\zzd$ so that $\lam\bigl(\<f\>\bigr) =\h(\af)$?
\end{question}

According to \cite{LSW}*{Lem.\ 4.3}, to compute the entropy of general algebraic $\zd$-actions it suffices to compute those defined by prime ideals. The case of principal prime ideals is part of Question \ref{ques:entropy}. By \cite{LSW}*{Thm.\ 4.2}, if $\mathfrak{p}$ is a nonprincipal prime ideal in $\zzd$ then $\h(\al_{\mathfrak{p}})=0$. 

\begin{question}\label{ques:nonprincipal}
   If $\mathfrak{p}$ is a nonprincipal prime ideal in $\zzd$, does $\lam(\mathfrak{p})=0$?
\end{question}

This approach to entropy suggests a formulation for algebraic actions of (not necessarily commutative) residually finite groups. Let $\Gamma$ be a countable discrete group. For each element $f$ in the integral group ring $\ZZ[\Gamma]$ there is an associated principal algebraic $\Gamma$-action $\af$ dual to the standard left-action of $\Gamma$ on $\ZZ[\Gamma]/\ZZ[\Gamma]f$. Assume now that $\Gamma$ is residually finite, and let $\mathcal{H}$ denote the collection of all finite-index subgroups of $\Gamma$. We say that $H\to\infty$ in $\mathcal{H}$ if for every finite subset $K\subset\Gamma$ we have eventually $H\cap K \subset\{1_{\Gamma}\}$. The quotient maps $\Gamma\to\Gamma/H$ for $H\in\mathcal{H}$ provide sofic approximations to $\Gamma$ for computing sofic entropy $\h(\af)$ of $\af$. There is an obvious extension of length $L$ to this situation.

We confine our attention to those $f$ which genuinely involve all of $\Gamma$. To do this, let $f\in\ZZ[\Gamma]$ have support $A=\supp f$. Say that $f$ has \emph{finite index} if $A^{-1}A := \{b^{-1}a:a,b\in A\}$ generates a finite-index subgroup of $\Gamma$. This property is invariant under left translation be elements of $\Gamma$.

\begin{question}
	With the previous notations, does
	\begin{displaymath}
	   \limsup_{H\to\infty} \frac{1}{|\Gamma/H|} \log L(\ZZ[\Gamma]f\cap\ZZ[H]) = \h(\af)
	\end{displaymath}
	for every finite-index $f\in \ZZ[\Gamma]$?
\end{question}

\appendix

\section{Computing the decimation limit of $1+x+y$}

There are few explicit calculations of the logarithmic Mahler measure, or more generally of the Ronkin function, of polynomials in $\zzd$ when $d\ge2$. Depending on the relative sizes of the coefficients, evaluation of the integrals involved typically requires the torus to be subdivided into a large number of subregions with complicated boundaries, and so simple formulas in terms of familiar functions are rare.

Here we treat the case $f(x,y)=1+x+y$ from Example \ref{exam:ledrappier}, where these calculations can be carried out, resulting in the formulas \eqref{eqn:led1} and \eqref{eqn:led2} for $\df$.

Smyth \cite{Smyth} first computed the logarithmic Mahler measure $\mahler(f)=\rf(0,0)$ to have the value in \eqref{eqn:smyth}. Twenty years later Maillot \cite{Mai}*{\S7.3}, aided by Cassigne, computed the entire Ronkin function $\rf(u,v)$, providing in his long memoir a concrete example of the canonical height of a hypersurface. Their result involves the Bloch-Wigner dilogarithm function, which is an alternative formulation of the series representation in our formulas. Lundqvist \cite{Lundqvist} gave the formulas for the partial derivatives of $\rf$ we use here. He also investigated the polynomial $1+x+y+z$, and showed that the second order partial derivatives of its Ronkin function can be expressed in terms of standard elliptic functions.

Let $\Del=\nf$ be the unit simplex, and denote its interior by $\Delo$. Let $\Af$ be the amoeba of $f$, as shown in Figure \ref{fig:amoeba}, and $\Afo$ be its interior. To evaluate $\rfs(r,s)$ for $(r,s)\in\Delo$, we need to know the value of $(u,v)\in\Afo$ at which the partial derivatives of $\rf(u,v)$ with respect to $u$ and $v$ equal $r$ and $s$, respectively. Fortunately, there is a simple relationship that was established by Lundqvist \cite{Lundqvist}, whose treatment we follow.

\begin{figure}[bt]
  \centering
  \begin{tikzpicture}[scale=1.5]
		\draw (0,0) -- (3,0) -- (2,2) -- cycle;
		\draw (2.7,1.2) node{\small $e^v$};
		\draw (0.8,1.2) node{\small $e^u$};
		\draw (1.5, 0) node[below]{\small 1};
		\draw (1.7, -0.6) node{(a)};
		\draw (2.65,0) node[above]{\small $\pi r$};
		\draw (.5,0) node[above]{\small $\pi s$};
		\draw (2.4,0) arc[start angle=180, end angle=117, radius = 0.6];
		\draw (.8,0) arc[start angle=0, end angle=40, radius=0.9];
  \end{tikzpicture}
  \hfil
  \begin{tikzpicture}[scale=3, >=Stealth]
  	  \draw [->] (0,-1.2) -- (0,1.2);
  	  \draw [->] (-.3, 0) -- (1.2,0);
  	  \draw[dashed]  (0,-1) arc[start angle=-90, end angle=90, radius=1];
  	  \draw  (0,0) -- (.65,0) -- (0.866,.5) -- cycle;
  	  \draw  (0,0) -- (.65,0) -- (0.866,-.5) -- cycle;
  	  \draw (.4,-1.3) node{(b)};
  	  \draw (.4,.35) node{\small $e^u$};
  	  \draw (.83,.2) node{\small $e^v$};
  	  \draw [->] (0.55,0) arc[start angle=180, end angle=63, radius=.1];
  	  \draw [->] (0.55,0) arc[start angle=180, end angle=297, radius=.1];
  	  \draw (.5, .12) node{\footnotesize $2\pi r$};
  	  \draw (.4,0) node[below]{\small 1};
    \end{tikzpicture}
    \caption{Determining partial derivatives from angles and sides}\label{fig:triangle}
\end{figure}
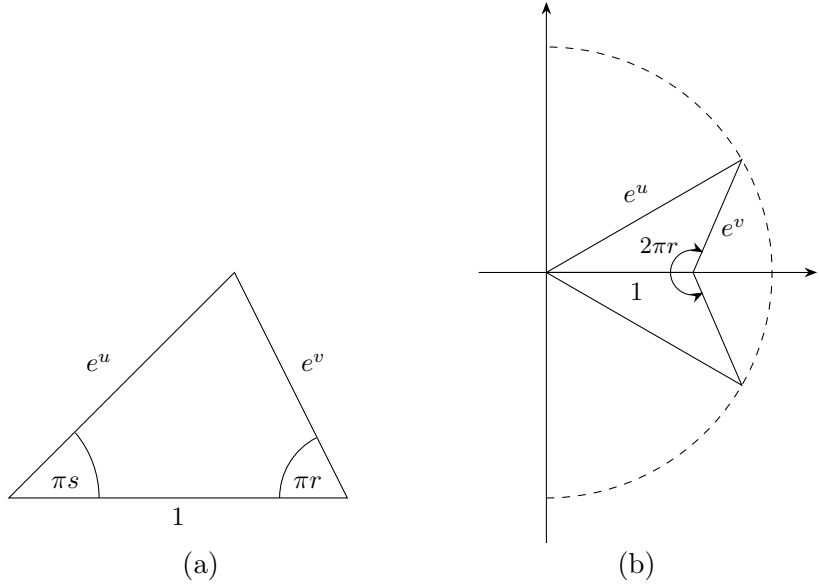

\begin{lemma}
	Let $(u,v)\in\Afo$, so that $1$, $e^u$, and $e^v$ form the sides of a nondegenerate triangle. Let $\pi r$ and $\pi s$ be the angles in this triangle shown in Figure \ref{fig:triangle}(a). Then
	\begin{equation}\label{eqn:partials}
	   \frac{\partial \rf}{\partial u}(u,v)= r \text{\quad and \quad} \frac{\partial \rf}{\partial v}(u,v)= s.
	\end{equation}
\end{lemma}

\begin{proof}
   We will compute the partial derivatives by differentiating the integrand in
   \begin{align*}
      \rf(u,v) &= \int_0^1 \int_0^1 \log |1+e^u \et + e^v\ep|\,d\theta\, d\phi\\
               &= \Re\Bigl[ \int_0^1 \int_0^1 \log (1+e^u \et + e^v\ep)\,d\theta\, d\phi\Bigr].
   \end{align*}
   In the last line $\log$ represents a local inverse to $\exp$, which is well-defined up to the addition of an integral multiple of $2 \pi i$. After taking partial derivatives, we will get a result that is independent of this multiple.
   
   By symmetry, it suffices to compute $\partial\rf / \partial u$. Differentiating the integrand gives
   \begin{displaymath}
      \frac{\partial \rf}{\partial u}(u,v) = \Re\Bigl[ \int_0^1\int_0^1
          \frac{e^u\et}{1+e^u\et + e^v\ep}\,d\theta\, d\phi \Bigr].
   \end{displaymath}
   Rewriting the integrals as contour integrals, we see that
   \begin{gather*}
      \int_0^1\int_0^1 \frac{e^u\et}{1+e^u\et +  e^v\ep}\,d\theta \,d\phi 
            = \frac{1}{(2 \pi i)^2} \int_{|z|=e^u} \int_{|w|=e^v} \frac{1}{1+z+w} \,dz \frac{dw}{w} \\
            = \frac{1}{2 \pi i} \int_{|w|=e^v} \Bigl[ \frac{1}{2 \pi i} \int_{|z|=e^u}
               \frac{dz}{z- (-1-w)} \Bigr]\,\frac{dw}{w}.
   \end{gather*}
   The inner integral is the winding number of the circle of radius $e^u$ around $-1-w=-1-e^v\ep$, and so has value $1$ if $|1+e^v\ep|<e^u$ and $0$ if $|1+e^v\ep|>e^u$ (these are mistakenly reversed in \cite{Lundqvist}). A glance at Figure \ref{fig:triangle}(b) shows that the value is $1$ for an interval of $\phi$ of length $2 \pi r$, and $0$ otherwise. Since 
   $ (1/2\pi i) (dw/w)$ is normalized Lebesgue measure ~$d\phi$, we obtain that $(\partial\rf/\partial u)(u,v) = r$.
\end{proof}

To compute the decimation limit $\df$, we need to express $u$ and $v$ in terms of $r$ and $s$. Let $a =e^u$ and $b=e^v$ be the sides of the triangle in Figure \ref{fig:triangle}(a). By the law of sines,
\begin{displaymath}
   \frac{a}{\sin \pi r} = \frac{b}{\sin \pi s} = \frac{1}{\sin \pi(1-r-s)}=\frac{1}{\sin \pi(r+s)},
\end{displaymath}
and hence
\begin{gather}
   a = a(r,s) = e^{u(r,s)}=\frac{\sin \pi r}{\sin\pi(r+s)}, \label{eqn:cova} \\
   b = b(r,s) = e^{v(r,s)} = \frac{\sin \pi s}{\sin\pi(r+s)}.\label{eqn:covb}
\end{gather}

For $(u,v)\in\Afo$ it follows from the definition \eqref{eqn:dual} that 
\begin{displaymath}
   -\rfs(u,v) = \inf_{(r,s)\in\Delo} \rf(u,v)-ru-sv,
\end{displaymath}
and by calculus the infimum is attained at the $(u,v)$ given by \eqref{eqn:partials}. Thus for $(r,s)\in\Delo$ we have that
\begin{equation}\label{eqn:}
   \df(r,s)=-\rfs\bigl( u(r,s),v(r,s)\bigr) = \rf\bigl(u(r,s),v(r.s)\bigr)-r\,u(r,s)
   - s\,v(r,s),
\end{equation}
where $u(r,s)$ and $v(r,s)$ are determined by \eqref{eqn:cova} and \eqref{eqn:covb}.

\begin{remark}
   Observe that the functions $u(r,s)$ and $v(r,s)$ in  \eqref{eqn:cova} and \eqref{eqn:covb} are real analytic on $\Delo$. Also, $\rf(u,v)$ is real analytic on $\Afo$. Together these show that $\df(r,s)$ is real analytic on $\Delo$.
\end{remark}

It remains to compute $\rf(u,v)$. By symmetry it suffices to assume that $u\ge v$. Using Jensen's formula \eqref{eqn:jensen}, we see that
\begin{align*}
   \rf(u,v) &= \int_0^1 \int_0^1 \log |1+e^u \et + e^v\ep|\,d\theta\, d\phi\\
            &= u+ \int_0^1 \int_0^1 \log |e^{-u}+e^{v-u}\ep + \et|\,d\theta \,d\phi\\
            &= u +\int_0^1 \log^+|e^{-u} + e^{v-u}\ep|\,d\phi.
\end{align*}
Note that $|e^{-u}+e^{v-u}\ep|\ge1$ if and only if $|1+e^v\ep|\ge e^u$, and another glance at Figure \ref{fig:triangle}(b) shows this occurs exactly when $-\pi(1-r)\le 2 \pi \phi\le \pi(1-r)$. Hence
\begin{align*}
   \rf(u,v) &= u + \int_{-\frac12(1-r)}^{\frac12(1-r)} \log |e^{-u}+e^{v-u}\ep|\,d\phi \\
   		   &= u - (1-r)u + \int_{-\frac12(1-r)}^{\frac12(1-r)} \log|1+e^v\ep|\,d\phi \\
   		   &= r\,u + \int_{-\frac12(1-r)}^{\frac12(1-r)} \log|1+e^v\ep|\,d\phi.
\end{align*}

First suppose that $e^v<1$, which corresponds to $(r,s)\in\Delo_1$, where $\Del_1^{}$ is defined in \eqref{eqn:delta1}. The series expansion of $\log(1+z)$ for $1+z$ in the domain of integration converges uniformly, and the imaginary part vanishes by symmetry. Hence
\begin{align*}
   \rf(u,v)&=r\,u + \int_{-\frac12(1-r)}^{\frac12(1-r)} \sum_{n=1}^\infty
       \frac{(-1)^{n+1}}{n} e^{nv} e^{2\pi i n \phi} \,d\phi\\
       &= r\,u + \sum_{n=1}^\infty \frac{(-1)^{n+1}}{n} e^{nv}\frac{1}{\pi n}
       \sin[\pi n(1-r)].
\end{align*}
Recalling that $e^{v(r,s)}=b(r,s)=(\sin \pi s)/\sin[\pi(r+s)]$, we conclude that
\begin{equation}\label{eqn:dfdelta1}
   \begin{split}
		\df(r,s)&=\rf\bigl(u(r,s),v(r,s)\bigr)-r\,u(r,s) - s\,v(r,s) \\
		        &= \sum_{n=1}^\infty \frac{(-1)^{n+1}}{\pi n^2} b(r,s)^n \sin[\pi n(1-r)]
		        - s\,\log[b(r,s)].
   \end{split}
\end{equation}

Now suppose that $e^v>1$, which corresponds to $(r,s)\in\Delo_2$, where $\Del_2$ is defined by \eqref{eqn:delta2}. Then $\log|1+e^v\ep|=v+\log|1+e^{-v}e^{-2\pi i\phi}|$. Calculating as before, 
\begin{align*}
   \rf(u,v)&=r\,u + (1-r)v + \int_{-\frac12(1-r)}^{\frac12(1-r)} \sum_{n=1}^\infty
       \frac{(-1)^{n+1}}{n} e^{-nv} e^{-2\pi i n \phi} \,d\phi\\
       &= r\,u + (1-r)v + \sum_{n=1}^\infty \frac{(-1)^{n+1}}{\pi n^2} b(r,s)^{-n}\sin[\pi n(1-r)].
\end{align*}
Thus for $(r,s)\in\Delo_2$ we find that
\begin{equation}\label{eqn:dfdelta2}
    \df(r,s) = \sum_{n=1}^\infty \frac{(-1)^{n+1}}{\pi n^2} b(r,s)^{-n} \sin[\pi n(1-r)]
    		        +(1-r-s)\log[b(r,s)].
\end{equation}
Finally, note that on the overlap $\Delta_1\cap\Delta_2$ inside $\Delo$, we have that $b(r,s)=1$ and so the series in \eqref{eqn:delta1} and \eqref{eqn:delta2} converge and agree, hence give the value of $\df(r,s)$ by continuity of the Legendre transform.


\begin{bibdiv}
\begin{biblist}

\bib{AtiyahMacdonald}{book}{
   author={Atiyah, M. F.},
   author={Macdonald, I. G.},
   title={Introduction to commutative algebra},
   publisher={Addison-Wesley Publishing Co., Reading, Mass.-London-Don
   Mills, Ont.},
   date={1969},
   pages={ix+128},
   review={\MR{0242802}},
}

\bib{Baker}{book}{
	author={Baker, Alan},
	title={Transcendental number theory},
	publisher={Cambridge Univ. Press, Cambridge},
	date={1975; reprinted 2022 with additional material},
	pages={xiv+176},
	review={\MR{0422171}}
}

\bib{Beer}{article}{
   author={Beer, Gerald},
   title={A natural topology for upper semicontinuous functions and a Baire
   category dual for convergence in measure},
   journal={Pacific J. Math.},
   volume={96},
   date={1981},
   number={2},
   pages={251--263},
   issn={0030-8730},
   review={\MR{637972}},
}
\bib{Berstel}{article}{
   author={Berstel, Jean},
   author={Mignotte, Maurice},
   title={Deux propri\'{e}t\'{e}s d\'{e}cidables des suites r\'{e}currentes lin\'{e}aires},
   language={French},
   journal={Bull. Soc. Math. France},
   volume={104},
   date={1976},
   number={2},
   pages={175--184},
   issn={0037-9484},
   review={\MR{414475}},
}
	

\bib{BoydUniform}{article}{
   author={Boyd, David W.},
   title={Uniform approximation to Mahler's measure in several variables},
   journal={Canad. Math. Bull.},
   volume={41},
   date={1998},
   number={1},
   pages={125--128},
   issn={0008-4395},
   review={\MR{1618904}},
   doi={10.4153/CMB-1998-019-6},
}



\bib{BoyleLind}{article}{
   author={Boyle, Mike},
   author={Lind, Douglas},
   title={Expansive subdynamics},
   journal={Trans. Amer. Math. Soc.},
   volume={349},
   date={1997},
   number={1},
   pages={55--102},
   issn={0002-9947},
   review={\MR{1355295}},
   doi={10.1090/S0002-9947-97-01634-6},
}
		

\bib{Burgos}{article}{
   author={Burgos Gil, Jos\'{e} Ignacio},
   author={Philippon, Patrice},
   author={Sombra, Mart\'{\i}n},
   title={Arithmetic geometry of toric varieties. Metrics, measures and
   heights},
   language={English, with English and French summaries},
   journal={Ast\'{e}risque},
   number={360},
   date={2014},
   pages={vi+222},
   issn={0303-1179},
   isbn={978-2-85629-783-4},
   review={\MR{3222615}},
}

\bib{Cassels}{book}{
   author={Cassels, J. W. S.},
   title={Local fields},
   publisher={Cambridge Univ. Press, Cambridge},
   date={1986},
   pages={xiv+360},
   review={\MR{861410}},
}



\bib{Chab}{article}{
   author={Chabauty, Claude},
   title={Limite d'ensembles et g\'{e}om\'{e}trie des nombres},
   language={French},
   journal={Bull. Soc. Math. France},
   volume={78},
   date={1950},
   pages={143--151},
   issn={0037-9484},
   review={\MR{38983}},
}
		

\bib{Conrad}{webpage}{
   author={Conrad, Keith},
   title={Galois descent},
   url={https://kconrad.math.uconn.edu/blurbs/galoistheory/galoisdescent.pdf},
}

\bib{Corn}{article}{
   author={Cornulier, Yves},
   title={On the Chabauty space of locally compact abelian groups},
   journal={Algebr. Geom. Topol.},
   volume={11},
   date={2011},
   number={4},
   pages={2007--2035},
   issn={1472-2747},
   review={\MR{2826931}},
   doi={10.2140/agt.2011.11.2007},
}
	

\bib{Duval}{article}{
   author={Duval, Dominique},
   title={Absolute factorization of polynomials: a geometric approach},
   journal={SIAM J. Comput.},
   volume={20},
   date={1991},
   number={1},
   pages={1--21},
   issn={0097-5397},
   review={\MR{1082133}},
   doi={10.1137/0220001},
}

\bib{EinsiedlerLind}{article}{
   author={Einsiedler, Manfred},
   author={Lind, Douglas},
   author={Miles, Richard},
   author={Ward, Thomas},
   title={Expansive subdynamics for algebraic ${\Bbb Z}^d$-actions},
   journal={Ergodic Theory Dynam. Systems},
   volume={21},
   date={2001},
   number={6},
   pages={1695--1729},
   issn={0143-3857},
   review={\MR{1869066}},
   doi={10.1017/S014338570100181X},
}
	

\bib{recurring}{book}{
   author={Everest, Graham},
   author={van der Poorten, Alf},
   author={Shparlinski, Igor},
   author={Ward, Thomas},
   title={Recurrence sequences},
   series={Mathematical Surveys and Monographs},
   volume={104},
   publisher={American Mathematical Society, Providence, RI},
   date={2003},
   pages={xiv+318},
   isbn={0-8218-3387-1},
   review={\MR{1990179}},
   doi={10.1090/surv/104},
}
		

\bib{Gao}{article}{
   author={Gao, Shuhong},
   title={Absolute irreducibility of polynomials via Newton polytopes},
   journal={J. Algebra},
   volume={237},
   date={2001},
   number={2},
   pages={501--520},
   issn={0021-8693},
   review={\MR{1816701}},
   doi={10.1006/jabr.2000.8586},
}
  
\bib{GKZ}{book}{
   author={Gel\cprime fand, I. M.},
   author={Kapranov, M. M.},
   author={Zelevinsky, A. V.},
   title={Discriminants, resultants, and multidimensional determinants},
   series={Mathematics: Theory \& Applications},
   publisher={Birkh\"{a}user Boston, Inc., Boston, MA},
   date={1994},
   pages={x+523},
   isbn={0-8176-3660-9},
   review={\MR{1264417}},
   doi={10.1007/978-0-8176-4771-1},
}



\bib{Gelfond}{book}{
   author={Gelfond, A. O.},
   title={Transcendental and algebraic numbers},
   publisher={Dover, New York},
   date={1960},
}

\bib{Gorin}{book}{
	author={Vadim Gorin},
	title={Lectures on random lozenge tilings},
	publisher={Cambridge Univ. Press},
	date={to appear},
}

\bib{Halmos}{article}{
   author={Halmos, Paul R.},
   title={On automorphisms of compact groups},
   journal={Bull. Amer. Math. Soc.},
   volume={49},
   date={1943},
   pages={619--624},
   issn={0002-9904},
   review={\MR{0008647}},
   doi={10.1090/S0002-9904-1943-07995-5},
}



\bib{HP}{article}{
   author={Hubbard, John},
   author={Pourezza, Ibrahim},
   title={The space of closed subgroups of ${\bf R}^{2}$},
   journal={Topology},
   volume={18},
   date={1979},
   number={2},
   pages={143--146},
   issn={0040-9383},
   review={\MR{544155}},
   doi={10.1016/0040-9383(79)90032-6},
}
		
\bib{KKS}{article}{
   author={Katok, Anatole},
   author={Katok, Svetlana},
   author={Schmidt, Klaus},
   title={Rigidity of measurable structure for ${\Bbb Z}^d$-actions by
   automorphisms of a torus},
   journal={Comment. Math. Helv.},
   volume={77},
   date={2002},
   number={4},
   pages={718--745},
   issn={0010-2571},
   review={\MR{1949111}},
   doi={10.1007/PL00012439},
}


\bib{KenyonOkounkovSheffield}{article}{
   author={Kenyon, Richard},
   author={Okounkov, Andrei},
   author={Sheffield, Scott},
   title={Dimers and amoebae},
   journal={Ann. of Math. (2)},
   volume={163},
   date={2006},
   number={3},
   pages={1019--1056},
   issn={0003-486X},
   review={\MR{2215138}},
   doi={10.4007/annals.2006.163.1019},
}



\bib{Ledrappier}{article}{
   author={Ledrappier, Fran\c{c}ois},
   title={Un champ markovien peut \^{e}tre d'entropie nulle et m\'{e}langeant},
   language={French, with English summary},
   journal={C. R. Acad. Sci. Paris S\'{e}r. A-B},
   volume={287},
   date={1978},
   number={7},
   pages={A561--A563},
   issn={0151-0509},
   review={\MR{512106}},
}
		

\bib{LSW}{article}{
   author={Lind, Douglas},
   author={Schmidt, Klaus},
   author={Ward, Tom},
   title={Mahler measure and entropy for commuting automorphisms of compact
   groups},
   journal={Invent. Math.},
   volume={101},
   date={1990},
   number={3},
   pages={593--629},
   issn={0020-9910},
   review={\MR{1062797}},}

\bib{Lundqvist}{article}{
   author={Lundqvist, Johannes},
   title={An explicit calculation of the Ronkin function},
   language={English, with English and French summaries},
   journal={Ann. Fac. Sci. Toulouse Math. (6)},
   volume={24},
   date={2015},
   number={2},
   pages={227--250},
   issn={0240-2963},
   review={\MR{3358612}},
   doi={10.5802/afst.1447},
}

\bib{MS}{book}{
   author={Maclagan, Diane},
   author={Sturmfels, Bernd},
   title={Introduction to tropical geometry},
   series={Graduate Studies in Mathematics},
   volume={161},
   publisher={American Mathematical Society, Providence, RI},
   date={2015},
   pages={xii+363},
   isbn={978-0-8218-5198-2},
   review={\MR{3287221}},
   doi={10.1090/gsm/161},
}

\bib{MagidinMcKinnon}{article}{
   author={Magidin, Arturo},
   author={McKinnon, David},
   title={Gauss's lemma for number fields},
   journal={Amer. Math. Monthly},
   volume={112},
   date={2005},
   number={5},
   pages={385--416},
   issn={0002-9890},
   review={\MR{2139573}},
   doi={10.2307/30037491},
}

\bib{MahlerCompactness}{article}{
   author={Mahler, K.},
   title={On lattice points in $n$-dimensional star bodies. I. Existence
   theorems},
   journal={Proc. Roy. Soc. London Ser. A},
   volume={187},
   date={1946},
   pages={151--187},
   issn={0962-8444},
   review={\MR{17753}},
   doi={10.1098/rspa.1946.0072},
}
	

\bib{Mahler}{article}{
   author={Mahler, K.},
   title={On some inequalities for polynomials in several variables},
   journal={J. London Math. Soc.},
   volume={37},
   date={1962},
   pages={341--344},
   issn={0024-6107},
   review={\MR{0138593}},
   doi={10.1112/jlms/s1-37.1.341},
}

\bib{Mai}{article}{
   author={Maillot, Vincent},
   title={G\'{e}om\'{e}trie d'Arakelov des vari\'{e}t\'{e}s toriques et fibr\'{e}s en droites
   int\'{e}grables},
   language={French, with English and French summaries},
   journal={M\'{e}m. Soc. Math. Fr. (N.S.)},
   number={80},
   date={2000},
   pages={vi+129},
   issn={0249-633X},
   review={\MR{1775582}},
   doi={10.24033/msmf.393},
}


\bib{Okounkov}{article}{
   author={Okounkov, Andrei},
   title={Limit shapes, real and imagined},
   journal={Bull. Amer. Math. Soc. (N.S.)},
   volume={53},
   date={2016},
   number={2},
   pages={187--216},
   issn={0273-0979},
   review={\MR{3474306}},
   doi={10.1090/bull/1512},
}

\bib{PassareRullgard}{article}{
   author={Passare, Mikael},
   author={Rullg\aa rd, Hans},
   title={Amoebas, Monge-Amp\`ere measures, and triangulations of the Newton
   polytope},
   journal={Duke Math. J.},
   volume={121},
   date={2004},
   number={3},
   pages={481--507},
   issn={0012-7094},
   review={\MR{2040284}},
   doi={10.1215/S0012-7094-04-12134-7},
}

\bib{Purbhoo}{article}{
   author={Purbhoo, Kevin},
   title={A Nullstellensatz for amoebas},
   journal={Duke Math. J.},
   volume={141},
   date={2008},
   number={3},
   pages={407--445},
   issn={0012-7094},
   review={\MR{2387427}},
   doi={10.1215/00127094-2007-001},
}

\bib{Rockafellar}{book}{
   author={Rockafellar, R. Tyrrell},
   title={Convex analysis},
   series={Princeton Mathematical Series, No. 28},
   publisher={Princeton University Press, Princeton, N.J.},
   date={1970},
   pages={xviii+451},
   review={\MR{0274683}},
}

\bib{Ronkin}{article}{
   author={Ronkin, L. I.},
   title={On zeros of almost periodic functions generated by functions
   holomorphic in a multicircular domain},
   language={Russian},
   conference={
      title={Complex analysis in modern mathematics (Russian)},
   },
   book={
      publisher={FAZIS, Moscow},
   },
   date={2001},
   pages={239--251},
   review={\MR{1833516}},
}

\bib{SchmidtBook}{book}{
   author={Schmidt, Klaus},
   title={Dynamical systems of algebraic origin},
   series={Progress in Mathematics},
   volume={128},
   publisher={Birkh\"auser Verlag, Basel},
   date={1995},
   pages={xviii+310},
   isbn={3-7643-5174-8},
   review={\MR{1345152}},
}

\bib{Smyth}{article}{
   author={Smyth, C. J.},
   title={On measures of polynomials in several variables},
   journal={Bull. Austral. Math. Soc.},
   volume={23},
   date={1981},
   number={1},
   pages={49--63},
   issn={0004-9727},
   review={\MR{615132}},
   doi={10.1017/S0004972700006894},
}

\bib{SmythKronecker}{article}{
   author={Smyth, C. J.},
   title={A Kronecker-type theorem for complex polynomials in several
   variables},
   journal={Canad. Math. Bull.},
   volume={24},
   date={1981},
   number={4},
   pages={447--452},
   issn={0008-4395},
   review={\MR{644534}},
   doi={10.4153/CMB-1981-068-8},
}


\end{biblist}
\end{bibdiv}
\end{document}